\documentclass[preprint]{elsarticle}

\usepackage{lineno}
\modulolinenumbers[5]
\usepackage{amsmath,amssymb,fixmath}
\usepackage{epsfig,epstopdf}
\usepackage{booktabs}
\usepackage{multirow}
\usepackage{subfig}
\usepackage{float}
\usepackage{mwe}
\usepackage{graphicx}
\usepackage{caption}
\usepackage{url}

\usepackage{array}
\usepackage{soul}

\usepackage{enumitem}
\usepackage{booktabs}
\usepackage{longtable}
\usepackage{stmaryrd}
\usepackage{multirow}
\usepackage{tabularx}
\usepackage{pdflscape}

\journal{Signal Processing - Journal - Elsevier}
\newcommand{\ds}{\displaystyle}

\begin{document}

\begin{frontmatter}

\title{Enhanced \emph{q}-Least Mean Square}


\author[mymainaddress]{Shujaat Khan}
\ead{shujaat@kaist.ac.kr}

\author[b]{Alishba Sadiq}
\ead{alishba.sadiq@pafkiet.edu.pk}

\author[b,c]{Imran Naseem\corref{mycorrespondingauthor}}
\cortext[mycorrespondingauthor]{Corresponding author}
\ead{imran.naseem@uwa.edu.au}

\author[c]{Roberto Togneri}
\ead{roberto.togneri@uwa.edu.au}

\author[d]{Mohammed Bennamoun}
\ead{mohammed.bennamoun@uwa.edu.au}

\address[mymainaddress]{Department of Bio and Brain Engineering, Korea Advanced Institute of Science and Technology (KAIST), Daejeon, Republic of Korea.}
\address[b]{College of Engineering, Karachi Institute of Economics and Technology,	Korangi Creek, Karachi 75190, Pakistan.}
\address[c]{School of Electrical, Electronic and Computer Engineering, The University of Western Australia, 35 Stirling Highway, Crawley, Western Australia 6009, Australia.}
\address[d]{School of Computer Science and Software Engineering, The University of Western Australia, 35 Stirling Highway, Crawley, Western Australia 6009, Australia.}

\begin{abstract}
In this work, a new class of stochastic gradient algorithm is developed based on $q$-calculus.  Unlike the existing $q$-LMS algorithm, the proposed approach fully utilizes the concept of $q$-calculus by incorporating time-varying $q$ parameter.  The proposed enhanced $q$-LMS ($Eq$-LMS) algorithm utilizes a novel, parameterless concept of error-correlation energy and normalization of signal to ensure high convergence, stability and low steady-state error.  The proposed algorithm automatically adapts the learning rate with respect to the error.  For the evaluation purpose the system identification problem is considered.  Extensive experiments show better performance of the proposed $Eq$-LMS algorithm compared to the standard $q$-LMS approach.

\end{abstract}

\begin{keyword}
Adaptive algorithms, Least Mean Squares Algorithm, $q$-calculus, Jackson derivative, system identification, $q$-LMS.
\end{keyword}

\end{frontmatter}

\linenumbers
\section{Introduction}\label{Sec:Intro}
Least square method is considered to be widely used optimization technique.  It has been applied in diversified applications such as plant identification \cite{sLMS}, detection of elastic inclusions \cite{W1}, noise cancellation \cite{NCLMS}, echo cancellation \cite{ECLMS}, ECG signal analysis \cite{ECGLMS}, elasticity imaging \cite{W3}, and time series prediction \cite{TSPLMS}, etc.  The least mean square (LMS) is one of the most popular least square algorithms for adaptive filtering due to its low computational complexity, however, it has a slow convergence rate due to the dependency on the eigenvalue-spread of the input correlation matrix \cite{NLMS}.  
Extensive research has been done towards the optimization of the LMS algorithm \cite{NCLMS,BcLMS,FOLMS,LMS1,qLMS}.  One of the disadvantages of the LMS is that it is sensitive to the scaling of its input.  In \cite{NLMS,LMS1}, the normalized LMS (NLMS) and its variants were proposed to solve this problem through normalization.  To improve the convergence rate and steady state performance, variable step size frameworks were devised in \cite{RVSS1,VSS}.  In \cite{CLMS,LMS9,LMS11,FCLMS}, different solutions for complex signal processing were proposed.  Similarly, to deal with non-linear signal processing problem, the concept of kernel function-based LMS algorithms was proposed in \cite{KLMS1,LMS6,LMS4}.  

Beside these variants, various definitions of gradient have also been used to derive improved LMS algorithms; for instance in \cite{RVSSFLMS}, a fractional order calculus (FOC) based least mean square algorithm, named as robust variable step size fractional least mean square (RVSS-FLMS), is proposed. The algorithm is derived using Riemann-Liouville fractional derivative for high convergence performance.  In \cite{VPFLMS,RVPFLMS}, some adaptive schemes were proposed for maintaining stability through adaptive variable fractional power.  The FOC variants are, however, not stable and diverge if the weights are negative or the input signal is complex \cite{bershad2017comments}.

Recently, the $q$-LMS algorithm is proposed which utilizes $q$-gradient from the Jackson's derivative so that the secant of the cost function is computed instead of the tangent \cite{qLMS}.  The algorithm takes larger steps towards the optimum solution and therefore, achieves a higher convergence rate.  The $q$-LMS algorithm has also been used for various applications, such as adaptive noise cancellation \cite{qLMS_ANC}, system identification, and designing of whitening filter \cite{qLMS_WF}.  In \cite{qNLMS} $q$-normalized LMS algorithm is proposed and its convergence performance is analyzed.  In \cite{qSSLMF}, using the same definition of $q$-calculus, variants of steady state least mean algorithms are derived.

All the aforementioned variants of the $q$-LMS algorithm enhance convergence speed at the cost of increased computational complexity and steady-state error.  In order to improve convergence rate without compromising the steady-state performance, a time-varying $q$-LMS is proposed in \cite{Tv_qLMS}.  However, it requires the tuning of two additional parameters ($\beta$ and $\gamma$), and the performance of the time-varying $q$-LMS \cite{Tv_qLMS} is very sensitive to the selection of the tuning parameters.  In this paper, we propose a new variant of the $q$-LMS by making the $q$-parameter time-varying.  The proposed enhanced $q$-LMS ($Eq$-LMS) utilizes a novel, parameterless concept of error-correlation energy and normalization to ensure rapid convergence without compromising stability and low steady-state error.  The proposed algorithm automatically adapts the learning rate with respect to error.  It takes larger steps in case of larger error and reduces the learning rate with decreased error.  Unlike the contemporary methods \cite{qLMS,RVSS1,VSS}, the proposed method is a parameterless technique and does not require manual tuning of any parameter.  The proposed algorithm is evaluated for the system identification problem and the results are demonstrated for both the steady-state performance and the convergence rate.  Extensive experiments are performed to show the superiority of the proposed $Eq$-LMS algorithm over a variety of contemporary methods.

The rest of the paper is structured as follows.  An overview of the $q$-calculus and $q$-least mean square algorithm is provided in Section \ref{overview_qcalculus}.  The details of proposed algorithm are discussed in Section \ref{Sec:Pro}, followed by the experimental results in Section \ref{Sec:Sim}.  The paper is finally concluded in Section \ref{Sec:Con}.

\section{Overview of \emph{q}-Least Mean Square Algorithm}\label{overview_qcalculus}
The conventional LMS algorithm is derived using the concept of steepest descent with the weight- update rule 
\begin{eqnarray}\label{LMS_WU}
\boldsymbol{w}(i+1) = \boldsymbol{w}(i) - \frac{\mu}{2} \nabla_{w}J(w),
\end{eqnarray}
where $J(w)$ is the cost function for the LMS algorithm and is defined as
\begin{equation}\label{cost_function}
J(w)  = E[e^{2}(i)],
\end{equation} 
where $E[\cdot]$ is the expectation operator and $e(i)$ is the estimation error between the desired response $d(i)$ and the output signal at the $i^{th}$ instant, i.e.,
\begin{eqnarray}
e(i) = d(i) - \boldsymbol{x}^{\intercal}(i) \boldsymbol{w}(i),
\end{eqnarray}
Here, $\boldsymbol{x}(i)$ is the input signal vector defined as
\begin{equation}
\boldsymbol{x}(i) = [{x}_{1}(i),{x}_{2}(i), \dots {x}_{M}(i)]^{\intercal},
\end{equation}
and $\boldsymbol{w}(i)$ is the weight vector defined as:
\begin{equation}
\boldsymbol{w}(i) = [w_{1}(i),w_{2}(i),\dots w_{M}(i)]
\end{equation}
where $M$ is the length of the filter.

\subsection{Overview of \emph{q}-Calculus}
The Quantum calculus or q-calculus is sometimes referred to as the calculus without a limit \cite{qbook1}.  It has been successfully used in various areas such as number theory, combinatorics, orthogonal polynomials, basic hyper-geometric functions and other sciences quantum theory, operational theory, mechanics, and the theory of relativity \cite{qbook2,qbook3,qbook4,qbook5}.

In q-calculus, the differential of a function is defined as (See, \cite{kac2001quantum})
\begin{equation}
d_{q}(f(x)) = f(qx)-f(x).
\end{equation}
The derivative therefore takes the form
\begin{equation}
D_{q}(f(x)) = \frac{d_{q}(f(x))}{d_{q}(x)} = \frac{f(qx)-f(x)}{(q-1)x}.
\end{equation}
When $q\rightarrow1$, the expression becomes the derivative in the classical sense.
The q-derivative of a function of the form $x^{n}$ is
\begin{equation}\label{q_derivative1}
D_{q,x}x^{n} =\begin{cases}
\ds\frac{q^{n}-1}{q-1} x^{n-1}, & q\neq1,
\\ 
\ds nx^{n-1}, & q = 1.
\end{cases}
\end{equation}
For a function $f(x)$ of $n$ number of variables, $\boldsymbol{x} = [ \boldsymbol{x}_{1}, \boldsymbol{x}_{2},.... \boldsymbol{x}_{n}]^{\intercal}$, the q-gradient is defined as
\begin{equation}\label{q_gradient1}
\nabla_{q,w} f(x)\triangleq [D_{q1,x1}f(x),D_{q2,x2}f(x),...D_{qn,xn}f(x)]^{\intercal} ,
\end{equation}
where $q =[q_{1},q_{2},\dots q_{N}]^{\intercal}$.

\subsection{\emph{q}-Least Mean Square (\emph{q}-LMS) Algorithm}
The performance of the LMS algorithm depends on the eigenvalue spread of the input correlation matrix.  The LMS is therefore regarded as an inherently slowly converging approach \cite{VSS}.  In order to resolve this issue the $q$-LMS has been proposed in \cite{qLMS}.  Instead of the conventional gradient, the $q$-LMS is derived using the $q$-calculus and utilizes the Jackson derivative method \cite{qLMS}, it takes larger steps (for $q>1$) in the search direction as it evaluates the secant of the cost function rather than the tangent \cite{qLMS}.  By replacing the conventional gradient in \eqref{LMS_WU} with the q-gradient, we get
\begin{eqnarray}\label{q_gradient}
\boldsymbol{w}(i+1) = \boldsymbol{w}(i) - \frac{\mu}{2} \nabla_{q,w}J(w).
\end{eqnarray}
The $q$-gradient of the cost function $J(w)$ for the $k^{th}$ weight is defined as
\begin{eqnarray}\label{q_derivative}
\nabla_{q,w_{k}}J(w) = \frac{\partial_{q_k}}{\partial_{q_k} e}J(w) \frac{\partial_{q_k}}{\partial_{q_k} y}e(i)\frac{\partial_{q_k}}{\partial_{q_k} w_{k}(i)}y(i).
\end{eqnarray}
Solving partial derivatives in \eqref{q_derivative} using the Jackson derivative defined in Section (\ref{overview_qcalculus}) gives
\begin{eqnarray}\label{first}
\frac{\partial_{q}}{\partial_{q} e}J(w) = \frac{\partial_{q}}{\partial_{q} e}(E[e^{2}(i)])= E[\frac{q_{k}^{2}-1}{q_{k}-1} e(i)]= E[(q_{k}+1) e(i)],
\end{eqnarray}
where $J(w)  = E[e^{2}(i)]$, $E[\cdot]$ is the expectation operator and $e(i) = d(i)-y(i)$.

Similarly
\begin{eqnarray}\label{second}
\frac{\partial_{q_k}}{\partial_{q_k} w_{k}(i)}y(i) = \boldsymbol{x}_{k}(i),
\end{eqnarray}
and
\begin{eqnarray}\label{third}
\frac{\partial_{q_k}}{\partial_{q_k} y}e(i) = -1,
\end{eqnarray}
Substituting equations \eqref{first}, \eqref{second}, and \eqref{third} in \eqref{q_derivative} gives
\begin{eqnarray}
\nabla_{q,w_{k}}(i)J(w) = -E[(q_{k}+1)e(i)\boldsymbol{x}_{k}(i)].
\end{eqnarray}
Similarly, for $k={1, 2, \dots, M}$,
\begin{equation}\label{q_derivative2}
\nabla_{q,w} J(w) = -E[(q_{1}+1)e(i)x_{1}(i),(q_{2}+1)e(i)\boldsymbol{x}_{2}(i), \dots (q_{M}+1)e(i)\boldsymbol{x}_{M}(i)].
\end{equation}

Consequently, Eq. \eqref{q_derivative2} can be written as
\begin{eqnarray}\label{q_gradient2}
\nabla_{q,w}J(w) = -2E[\mathbf{G} \boldsymbol{x}(i) e(i)],
\end{eqnarray}\label{diagonal_G}
where $\mu$ is the learning rate (step-size) and $\mathbf{G}$ is a diagonal matrix
\begin{eqnarray}\label{G}
{\rm diag}(\mathbf{G}) = [(\frac{q_{1}+1}{2}), (\frac{q_{2}+1}{2}),.....(\frac{q_{M}+1}{2})]^{\intercal}.
\end{eqnarray}

Dropping the expectation from the q-gradient in \eqref{q_gradient2} results in
\begin{eqnarray}\label{ins_val}
\nabla_{q,w}J(w) \approx -2G \boldsymbol{x}(i) e(i).
\end{eqnarray}
Substituting \eqref{ins_val} in \eqref{LMS_WU} renders the weight update rule of the q-LMS algorithm by
\begin{eqnarray}\label{qLMS_final}
\boldsymbol{w}(i+1) = \boldsymbol{w}(i) +\mu G \boldsymbol{x}(i)e(i).
\end{eqnarray}

\subsubsection{\emph{q}-LMS as a whitening filter-\emph{q}-Normalized LMS}
The $q$-normalized least mean square ($q$-NLMS) algorithm is defined (See \cite{qNLMS})
\begin{eqnarray}\label{qNLMS_final}
\boldsymbol{w}(i+1) = \boldsymbol{w}(i) +\mu \frac{G\boldsymbol{x}(i)e(i)}{\zeta + ||\boldsymbol{x}(i)||^2_G},
\end{eqnarray}
where $\zeta$ is a small value added in the denominator to avoid the indeterminate form, and $||\boldsymbol{x}(i)||^2_G$ is the weighted norm of the input vector.  By selecting the $q$ parameter in \eqref{qLMS_final} as $q = 1/\lambda\max$, we can design a whitening filter and hence it can remove the dependency on the input correlation \cite{qNLMS}.

\subsubsection{Time-varying \emph{q}-LMS}
The time-varying \emph{q}-LMS algorithm is based on the variable step size (VSS) method \cite{VSS} and is given (See \cite{Tv_qLMS})  
\begin{equation}
\Psi(i+1) = \beta\Psi(i) + \gamma e(i)^{2}, \qquad  (0<\beta<1, \gamma>0),
\end{equation}
\begin{equation}
q(i+1) = \begin{cases}
\ds q_{upper}, & \Psi(i+1)>q_{upper},
\\ 
\ds 1, & \Psi(i+1)<1,
\\
\ds \Psi & otherwise.
\end{cases}
\end{equation}

where $q_{upper}$ is so chosen to satisfy the stability bound, (See \cite{qLMS})
\begin{equation}
q_{upper} = \frac{2}{ \mu \lambda_{max}}.
\end{equation}

The $q(i+1)$ is updated according to the estimation of the square of the estimation error. When the estimation error is large, $q(i)$ will approach its upper bound denoted by $q_{upper}$, while for smaller values $q(i)$ goes to unity for a lower steady-state error.

\section{The Proposed Enhanced \emph{q}-Least Mean Square (Eq-LMS) Algorithm}\label{Sec:Pro}
The $q$-LMS algorithm has an extra degree of freedom to control the performance via the diagonal matrix G, which comprises of the $q$-dependent entries.  The weight-update rule of the $q$-LMS algorithm can be written as
\begin{eqnarray}\label{w_update}
\boldsymbol{w}(i+1) = \boldsymbol{w}(i) +\mu  \boldsymbol{\hat{x}}(i)e(i),
\end{eqnarray}
where $\boldsymbol{\hat{x}}(i) = G \boldsymbol{x}(i)$.  For the special case of G=I (identity matrix), the $q$-LMS algorithm will be transformed into the conventional LMS.  
Based on the above discussion, we make the following important observations.  

\begin{itemize}
	\item We argue that the $q$-gradient with $q>1$ enhances the speed of convergence as it takes the secant of function rather than the tangent \cite{qLMS}.  The larger the value of the $q$ parameter, the faster the convergence of the algorithm.  But this improvement in the rate of convergence comes at the cost of a degradation in the steady-state performance. 
	\item The time varying $q$-LMS technique \cite{Tv_qLMS} is based on variable step-size method \cite{VSS}, which requires the tuning of additional parameters such as $\beta$ and $\gamma$.
	\item[$\bullet$] By selecting the $q$ parameter in \eqref{qLMS_final} as $q = 1/\lambda\max$, we can design a whitening filter and hence it can remove the dependency on the input correlation \cite{qNLMS}.  However, with a large step-size the $q$-NLMS converges rapidly with a compromised steady state performance.  Similarly, a smaller step-size results in better steady state performance but with slow convergence.  As such, the two important performance parameters cannot be optimized simultaneously.	
\end{itemize}

\subsection{Proposed Improvements}
To overcome the aforementioned issues, we propose the $Eq$-LMS algorithm with the following improvements.
\begin{itemize}
  \item[$\bullet$] To achieve higher convergence rate with lower steady-state error, we propose to incorporate the instantaneous error energy to adapt the $q$-parameter.  The proposed algorithm automatically takes large steps when the error is large and reduces the step-size with the decreasing error.  Note that, unlike the time-varying $q$-LMS \cite{Tv_qLMS}, no additional tuning parameters are introduced and the proposed approach is completely automatic.
  
  \item[$\bullet$] The whitening factor ($q = 1/\lambda\max$) is also utilized to set the limits of adaptive $q$-vector, this allows the algorithm to operate at a higher convergence rate without worrying about the divergence issues.
  
  \item[$\bullet$] To update each $q$-parameter, the proposed $Eq$-LMS utilize, a \emph{responsible error} for each tap of the filter.  With this improvement, the $q$-variable for each tap will be updated accordingly, hence, both the steady-state error and convergence performance can be improved significantly.

\end{itemize}

\subsection{Formulation of the Eq-LMS Algorithm}
By replacing the fixed $\mathbf{G}$ in \eqref{qLMS_final} with its time-varying form $\mathbf{G}(i)$, the weight update rule of the proposed $Eq$-LMS is given as
\begin{equation}
\boldsymbol{w}(i+1) = \boldsymbol{w}(i) + \mu e(i) \boldsymbol{x}(i)\mathbf{G}(i),
\end{equation}
where $\mu$ is the learning rate, and $e(i)$ is the error at the $i^{th}$ instant defined by
\begin{equation}
e(i) = d(i) - y(i),
\end{equation}
where $d(i)$ and $y(i)$ are the desired and estimated output at the $i^{th}$ instant, respectively.
Here, $\mathbf{G}(i)$ is a diagonal matrix with time-varying diagonal elements, and is defined as
\begin{equation} 
\mathbf{G}(i) =
\begin{bmatrix}
q_{1}(i) & & \\
& q_{2}(i) & \\
& &	\ddots & \\
& & & q_{M}(i)
\end{bmatrix}.
\end{equation}
It can also be written as
\begin{equation}
	\mathbf{G}(i)={\rm diag}(q_1(i), \dots, q_M(i)) = {\rm diag}(\boldsymbol{q}(i)),
\end{equation}
where,
\begin{equation}
\boldsymbol{q} (i) = \{q_{1}(i),q_{2}(i), \dots q_{M}(i)\}.
\end{equation}
We propose an update rule for vector $\boldsymbol{q}$ defined in the following steps.
\begin{itemize}
	\item Step1: Initialize vector $\boldsymbol{q}$ with any positive random values.
	\item Step2: Use instantaneous error to update first entry $q_1$ of $\mathbf{q}$ vector, which is associated with weight of the instant input tap, i.e., 
	\begin{equation}
	q_{1}(i+1) = \frac{1}{M+1}\{|e(i)|+\sum_{j=1}^{M}q_{j}(i)\},
	\end{equation}
	where $M$ is the length of the filter.
	\item Step3: To avoid divergence while maintaining the higher convergence rate, the following conditions will be evaluated:
	\begin{equation}\label{q_parameter}
	\boldsymbol{q}(i+1) = \left\{ \begin{array}{rcl}
	\frac {1} {\lambda_{max}} & if & |q_{1}(i+1)| > \frac {1} {\lambda_{max}},\\ 
	q_{1}(i+1) & if & |q_{1}(i+1)| < \frac {1} {\lambda_{max}},\\ 
	\end{array}\right.
	\end{equation}
	where $\lambda_{max}$ is the maximum eigenvalue of the input auto-correlation matrix.  
	\item Step4: Update all entries of vector $\boldsymbol{\alpha}$ except for the first entry, simply by shifting:    
	\begin{equation}
	q_{k+1}(i+1) = q_{k}(i),
	\end{equation}
	where $1 < k < M-1$
	\item Step5: For next iterations, repeat steps 2 to 5.
\end{itemize}

Finally, the weight-update equation of the proposed $Eq$-LMS can be written as:
\begin{equation}\label{EqLMS_final}
\boldsymbol{w}(i+1) = \boldsymbol{x}(i) + \mu e(i)\boldsymbol{x}(i)\odot\boldsymbol{q}(i),
\end{equation}
where $\odot$ indicates the element wise multiplication.

\section{Experiments}\label{Sec:Sim}
\begin{figure}[h]
	\begin{center}
		\centering
		\includegraphics*[scale=0.4,bb=0 0 520 360]{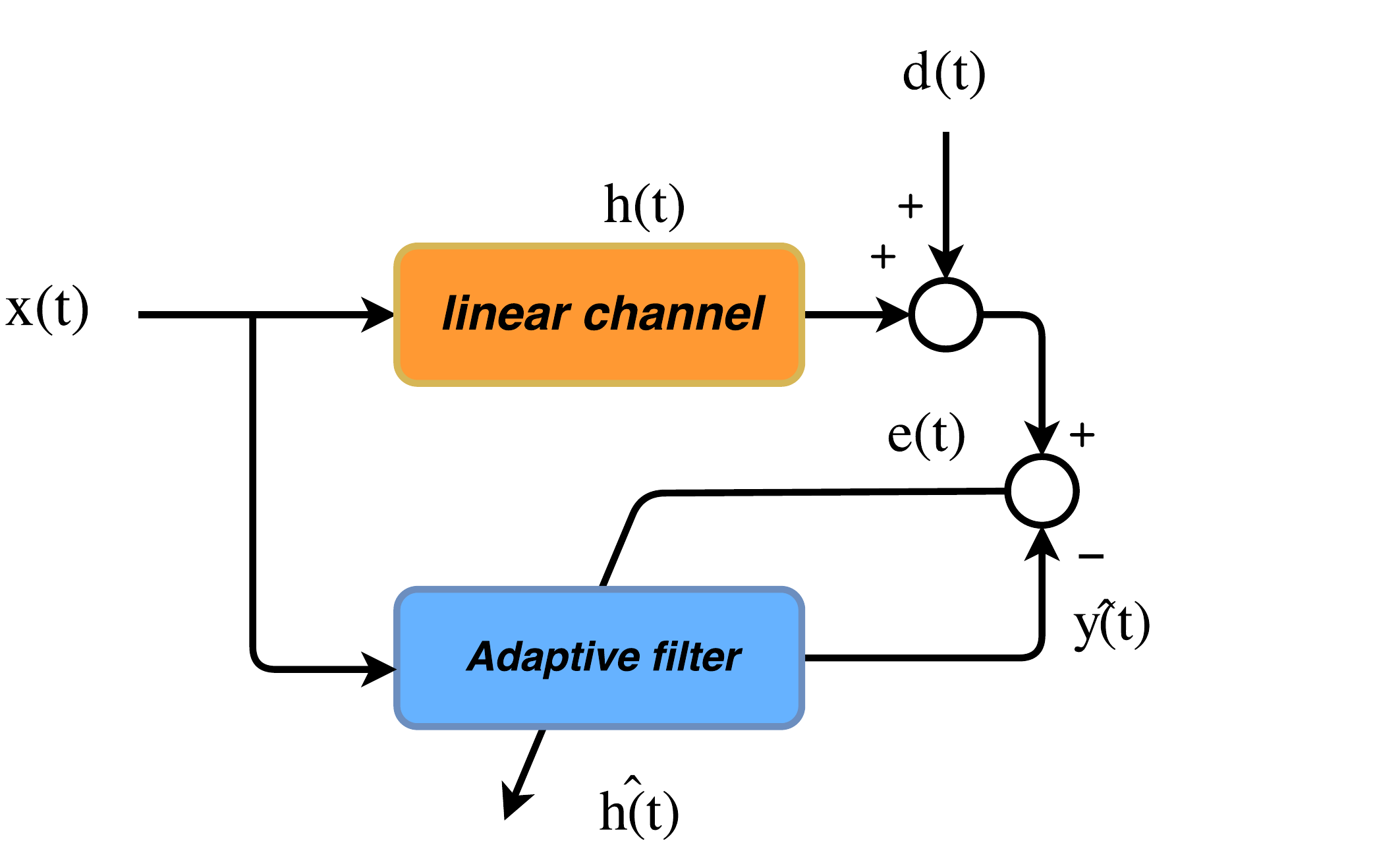}
	\end{center}
	\caption{Channel estimation using adaptive learning algorithm.}
	\label{plant1}
\end{figure}

For the evaluation of the proposed algorithm the problem of system identification is used.  Adaptive learning methods have been successfully used to identify the unknown system,  with numerous applications, for example, in control engineering, communication systems \cite{LMS1,FLMF, opfilter, khan2016novel, FBPTT,LMS14}.  
Channel estimation, for instance, is a widely used method in communication systems to estimate the characteristics of an unknown channel.  Consider a linear channel shown in Fig. \ref{plant1}.
\begin{equation}\label{math_model}
y(t) = h_{1}x(t)+h_{2}x(t-1)+h_{3}x(t-2)+h_{4}x(t-3)+h_{5}x(t-4).
\end{equation}
Equation \eqref{math_model} shows the mathematical model of the system, where $x(t)$ and $y(t)$ are the input and output of the system, respectively and $d(t)$ is the disturbance which is taken to be white Gaussian noise in this case.  For this experiment, $\boldsymbol{x}(t)$ is chosen to be consisting of $1\times10^{6}$ randomly generated samples obtained from Gaussian distribution of mean zero and variance of $1$.  In \eqref{math_model}, the system is defined by its impulse response $h(t)$ while $\hat{y}(t)$, $\hat{h}(t)$, and $e(t)$ are the estimated output, estimated impulse response, and the error of estimation, respectively.  The simulation parameters selected are as follows:  coefficient values of $h_1=-2$, $h_2=-1$, $h_3=0$, $h_4=1$ and $h_5=2$ are selected for the channel, the experiments are performed on three noise levels with the SNR values of $10$dB, $20$dB and $30$dB.  The weights are initialized to zero for all algorithms.  Specifically, the objective of these simulations is to compare the performance of the proposed enhanced $q-$LMS ($Eq$-LMS) algorithm with the contemporary counterparts, i.e., Least Mean Square (LMS)/ $q$-LMS \cite{qLMS} at $q=1$, $q$-LMS \cite{qLMS} at $q=2$, time-varying $q$-LMS \cite{Tv_qLMS} and the normalized LMS (NLMS), for given convergence rate and given steady state error in three different scenarios.

For the performance evaluation, the normalized weight deviation (NWD) in the actual and the obtained weights is compared. Specifically, we define
\begin{equation}\label{NWD}
\mbox{NWD}=\frac{\left\Vert \bf h-\bf w \right\Vert}{\left\Vert\bf h\right\Vert},
\end{equation} 
where $\boldsymbol{h}$ is the actual impulse response of the channel and $\boldsymbol{w}$ is the estimated weight-vector.  The simulations are repeated for $1000$ independent runs and mean results are reported.  The simulations are performed primarily to evaluate the steady-state and convergence performances of the proposed algorithm for various learning rates. Accordingly, three Evaluation protocols are designed: 

\begin{enumerate}
	\item {\bf Evaluation protocol 1}: learning rate=$1\times10^{-1}$, SNR=\{10,20,30\} dB.
	\item {\bf Evaluation protocol 2}: learning rate=$1\times10^{-2}$, SNR=\{10,20,30\} dB.
	\item {\bf Evaluation protocol 3}: learning rate=$1\times10^{-3}$, SNR=\{10,20,30\} dB.
\end{enumerate}

\subsection{Experiments to Evaluate the Steady-State Performance}
\subsubsection{Evaluation Protocol 1: Fast Convergence}
For the comparison of the steady-state performance,  all algorithms were setup for equal convergence and after $10000$ iterations, the steady-state value of the $\mbox{NWD}$ is examined.  The learning rate (step size) configurations for equal convergence rate (Evaluation Protocol 1) is shown in Table \ref{con_tab_C_scenario1}.
\begin{table}[H]
	\centering
	\caption{Evaluation protocol 1: Configuration of learning rates of different approaches for an equal convergence rate of $1 \times 10^{-1}$. }
\label{con_tab_C_scenario1}
	\begin{tabular}{|l|c|c|c|}
		\hline
		\multicolumn{1}{|c|}{\multirow{2}{*}{\textbf{Algorithm}}} & \multicolumn{3}{c|}{\textbf{Learning Rate $\mu$}} \\ \cline{2-4} 
		\multicolumn{1}{|c|}{} & \textbf{10 dB SNR} & \textbf{20 dB SNR} & \textbf{30 dB SNR} \\ \hline
		\textbf{LMS/$q$-LMS} & $4.4\times10^{-2}$ & $2.45\times10^{-2}$ & $2.7\times10^{-5}$ \\ \hline
		\textbf{Time-varying $q$-LMS} & $3.5\times10^{-1}$ & $1.3\times10^{-1}$ & $3.2\times10^{-5}$ \\ \hline
		\textbf{Normalized LMS} & $2.45\times10^{-1}$ & $1.2\times10^{-1}$ & $8.7\times10^{-5}$ \\ \hline
		\textbf{Proposed $Eq$-LMS} & $1\times10^{-1}$ & $1\times10^{-1}$ & $1\times10^{-1}$ \\ \hline
	\end{tabular}
\end{table}

The relevant normalized weight difference (NWD) curves with three different SNR values are depicted in Fig. \ref{convergence_scenario}.  From the Fig. \ref{convergence_scenario} (a),(b), and (c), it can be seen that the proposed $Eq$-LMS produced the best performance under all three conditions:  (1) for the SNR value of $10$dB, it outperformed the LMS/$q$-LMS at $q=1$, the $q$-LMS at $q=2$, the time-varying $q$-LMS, and the NLMS by the NWD value of $-1.03$ dB, $-2.95$ dB, $-8.06$ dB, and $-2.14$ dB, respectively, (2) for the SNR value of $20$ dB, it surpassed the listed algorithms by the NWD value of $-2.36$ dB, $-4.05$ dB, $-6.96$ dB, and $-3.29$ dB, respectively, and (3) for the SNR value of $30$ dB, the above mentioned algorithms were outperformed by a margin of $-2.29$ dB, $-3.86$ dB, $-7.18$, and $-3.47$ dB, respectively.  Note that the proposed $Eq$-LMS algorithm showed the lowest steady state error in all conditions while the LMS, NLMS and the $q$-LMS at $q=2$ show faster convergence than the time-varying $q$-LMS but with greater steady state error than the proposed $Eq$-LMS.  With the above discussed settings, results for the channel estimation problem are summarized in Table \ref{Res_tab_R_scenario1}. 
\begin{table}[H]
	\centering
	\caption{Evaluation protocol 1: Results of various approaches for an equal convergence rate.}
	\label{Res_tab_R_scenario1}
	\begin{tabular}{|l|c|c|c|}
		\hline
		\multicolumn{1}{|c|}{\multirow{2}{*}{\textbf{Algorithm}}} & \multicolumn{3}{c|}{\textbf{Steady-state NWD (dB)}} \\ \cline{2-4} 
		\multicolumn{1}{|c|}{} & \textbf{10 dB SNR} & \textbf{20 dB SNR} & \textbf{30 dB SNR} \\ \hline
		\textbf{LMS/$q$-LMS ($q=1$)} & -14.63 & -21.06 & -28.77 \\ \hline
		\textbf{$q$-LMS ($q=2$)} & -12.71 & -19.37 & -27.02 \\ \hline	
		\textbf{Time-varying $q$-LMS} & -7.60 & -16.46 & -23.88 \\ \hline
		\textbf{Normalized LMS} & -13.52 & -20.13 & -27.59 \\ \hline
		\textbf{Proposed $Eq$-LMS} & \textbf{-15.66} & \textbf{-23.42} & \textbf{-31.06} \\ \hline
	\end{tabular}
\end{table}

\subsubsection{Evaluation protocol 2: Medium Convergence}
The learning rate (step size) configurations for an equal convergence rate (Evaluation protocol 2) is shown in Table \ref{con_tab_C_scenario2}. 

\begin{table}[H]
	\centering
	\caption{Evaluation protocol 2: Configuration of learning rates of different approaches for an equal convergence rate of $1 \times 10^{-2}$. }
	\label{con_tab_C_scenario2}
	\begin{tabular}{|l|c|c|c|}
		\hline
		\multicolumn{1}{|c|}{\multirow{2}{*}{\textbf{Algorithm}}} & \multicolumn{3}{c|}{\textbf{Learning Rate $\mu$}} \\ \cline{2-4} 
		\multicolumn{1}{|c|}{} & \textbf{10 dB SNR} & \textbf{20 dB SNR} & \textbf{30 dB SNR} \\ \hline
		\textbf{LMS/$q$-LMS} & $4\times10^{-3}$ & $1.5\times10^{-3}$ & $5.2\times10^{-4}$ \\ \hline
		\textbf{Time-varying $q$-LMS} & $2\times10^{-2}$ & $1\times10^{-3}$ & $1.8\times10^{-3}$ \\ \hline
		\textbf{Normalized LMS} & $2\times10^{-2}$ & $9\times10^{-3}$ & $2.5\times10^{-3}$ \\ \hline
		\textbf{Proposed $Eq$-LMS} & $1\times10^{-2}$ & $1\times10^{-2}$ & $1\times10^{-2}$ \\ \hline
	\end{tabular}
\end{table}

The relevant normalized weight difference (NWD) curves with three different SNR values are depicted in Fig. \ref{convergence_scenario}.  From the Fig. \ref{convergence_scenario} (d),(e), and (f), it can be seen that the proposed $Eq$-LMS produced the best performance under all three conditions:  (1) for the SNR value of $10$ dB, it outperformed the LMS/$q$-LMS at $q=1$, $q$-LMS at $q=2$, time-varying $q$-LMS and the NLMS by the NWD value of $-0.68$ dB, $-2.23$ dB, $-4.08$ dB, and $-1.74$ dB, respectively, (2) for the SNR value of $20$ dB, it surpassed the above-mentioned algorithms by the NWD value of $-1.22$ dB, $-2.74$ dB, $-5.34$ dB, and $-2.45$ dB, respectively, and (3) for the SNR value of $30$ dB, the above mentioned algorithms were outperformed by a margin of $-1.44$ dB, $-2.94$ dB, $-4.07$ and $-2.47$ dB, respectively. Note that the proposed $Eq$-LMS algorithm showed the lowest steady state error in all conditions while the LMS, NLMS and the $q$-LMS at $q=2$ show faster convergence than time-varying $q$-LMS but with a greater steady state error than the proposed $Eq$-LMS.  With the above discussed settings, results for the channel estimation problem are summarized in Table \ref{Res_tab_R_scenario2}.   

\begin{table}[H]
	\centering
	\caption{Evaluation protocol 2: Results of various approaches for an equal convergence rate.}
	\label{Res_tab_R_scenario2}
	\begin{tabular}{|l|c|c|c|}
		\hline
		\multicolumn{1}{|c|}{\multirow{2}{*}{\textbf{Algorithm}}} & \multicolumn{3}{c|}{\textbf{Steady-state NWD (dB)}} \\ \cline{2-4} 
		\multicolumn{1}{|c|}{} & \textbf{10 dB SNR} & \textbf{20 dB SNR} & \textbf{30 dB SNR} \\ \hline
		\textbf{LMS/qLMS (q=1)} & -25.42 & -32.70 & -40.00 \\ \hline
		\textbf{qLMS (q=2)} & -23.62 & -30.86 & -38.19 \\ \hline
		\textbf{Time-varying qLMS} & -22.40 & -29.40 & -37.45 \\ \hline
		\textbf{Normalized LMS} & -24.18 & -31.45 & -38.80 \\ \hline
		\textbf{Proposed Eq-LMS} & \textbf{-25.85} & \textbf{-33.60} & \textbf{-41.13} \\ \hline
	\end{tabular}
\end{table}

\subsubsection{Evaluation protocol 3: Slow Convergence}
The learning rate (step size) configurations for equal convergence rate (Evaluation protocol 3) is shown in Table \ref{con_tab_C_scenario3}.
\begin{table}[H]
	\centering
\caption{Evaluation protocol 3: Configuration of learning rates of different approaches for the fix convergence rate of $1 \times 10^{-3}$. }
	\label{con_tab_C_scenario3}
	\begin{tabular}{|l|c|c|c|}
		\hline
		\multicolumn{1}{|c|}{\multirow{2}{*}{\textbf{Algorithm}}} & \multicolumn{3}{c|}{\textbf{Learning Rate $\mu$}} \\ \cline{2-4} 
		\multicolumn{1}{|c|}{} & \textbf{10 dB SNR} & \textbf{20 dB SNR} & \textbf{30 dB SNR} \\ \hline
		\textbf{LMS/qLMS} & $3.6\times10^{-4}$ & $1.3\times10^{-4}$ & $4.5\times10^{-3}$ \\ \hline
		\textbf{Time-varying qLMS} & $1.6\times10^{-3}$ & $6\times10^{-4}$ & $1.5\times10^{-4}$ \\ \hline
		\textbf{Normalized LMS} & $1.9\times10^{-3}$ & $7\times10^{-4}$ & $2.3\times10^{-4}$ \\ \hline
		\textbf{Proposed Eq-LMS} & $1\times10^{-3}$ & $1\times10^{-3}$ & $1\times10^{-3}$ \\ \hline
	\end{tabular}
\end{table}

\begin{figure}[]
	\begin{tabular}{>{\centering}m{3.7cm} >{\centering}m{3.7cm} >{\centering\arraybackslash}m{3.7cm}}

		\raisebox{-\totalheight}{\centering \fbox{\includegraphics*[scale=0.27,bb=100 240 500 550]{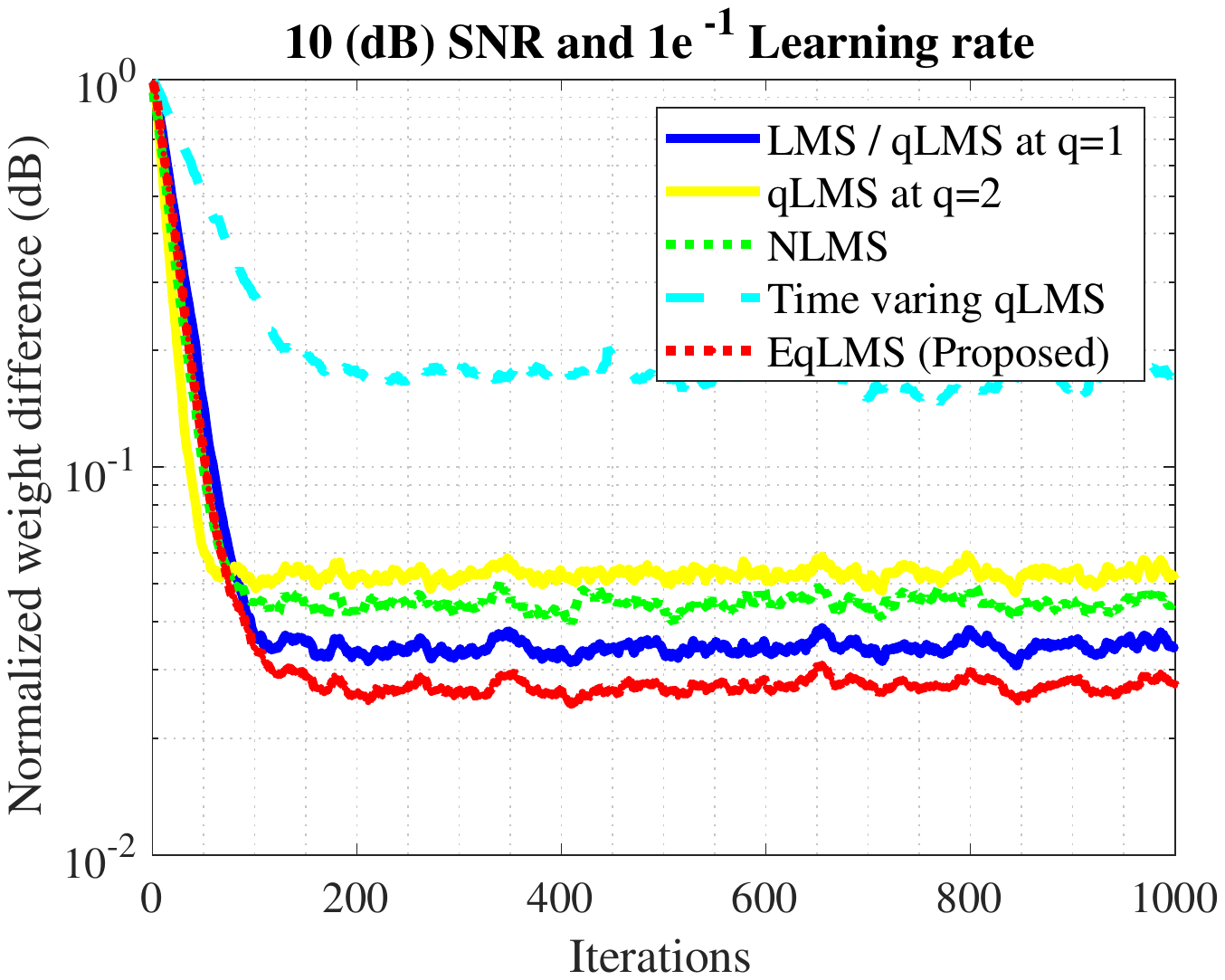}}}
		& 

		\raisebox{-\totalheight}{\centering \fbox{\includegraphics*[scale=0.27,bb=100 240 500 550]{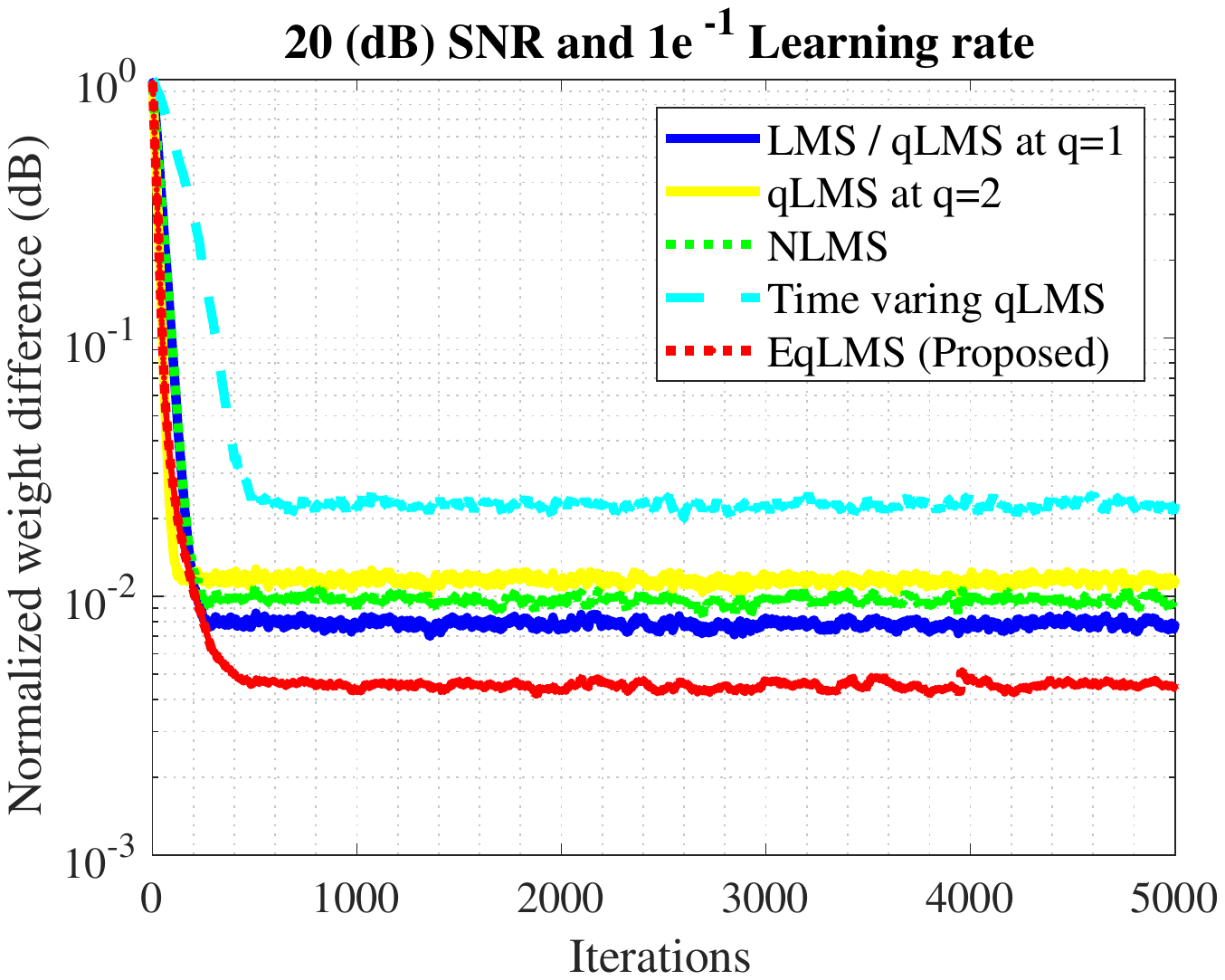}}}
		& 

		\raisebox{-\totalheight}{\centering \fbox{
				\includegraphics*[scale=0.27,bb=100 240 500 550]{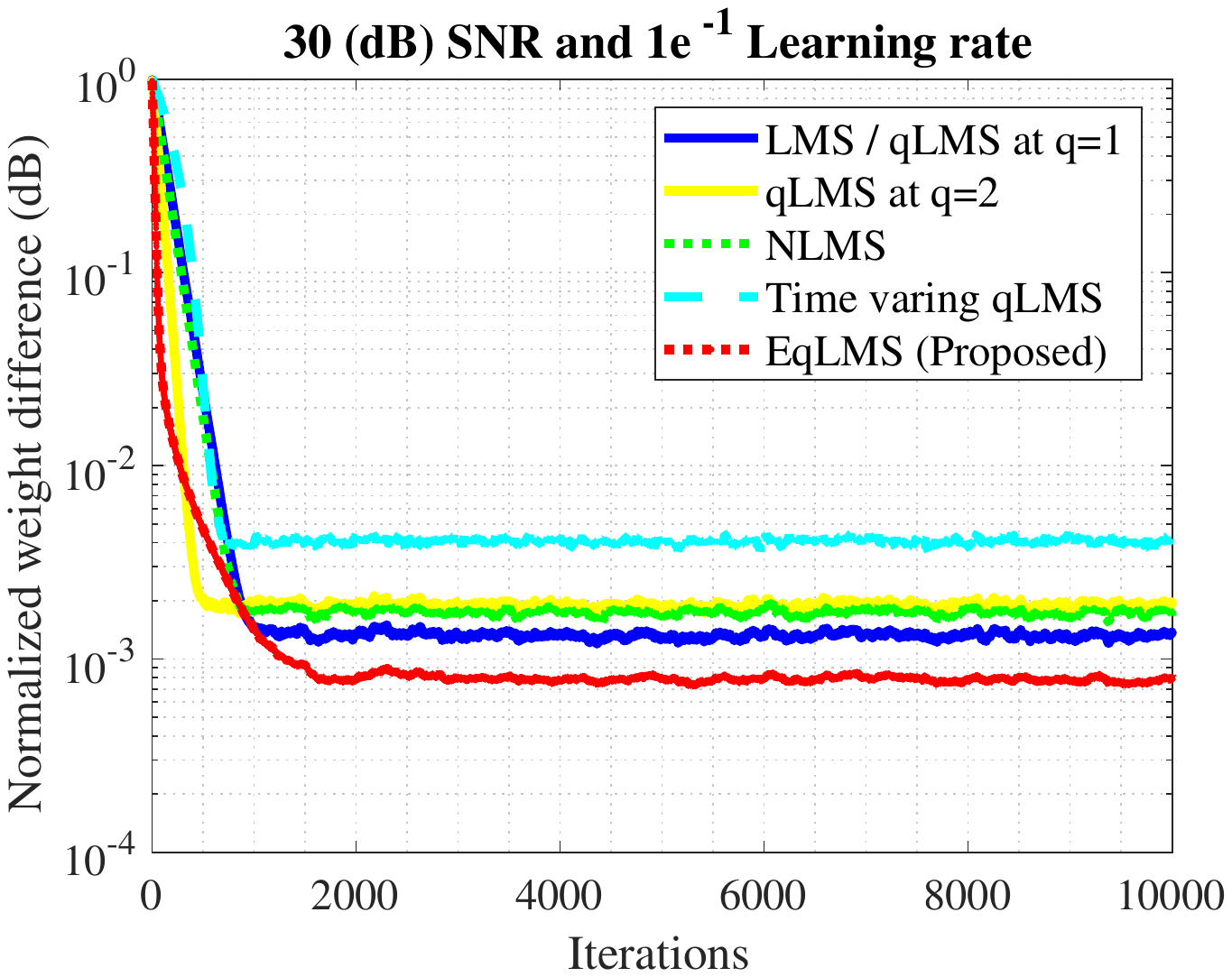}}}
		\\ 
		a & b & c \\

		\raisebox{-\totalheight}{\centering \fbox{\includegraphics*[scale=0.27,bb=100 240 500 550]{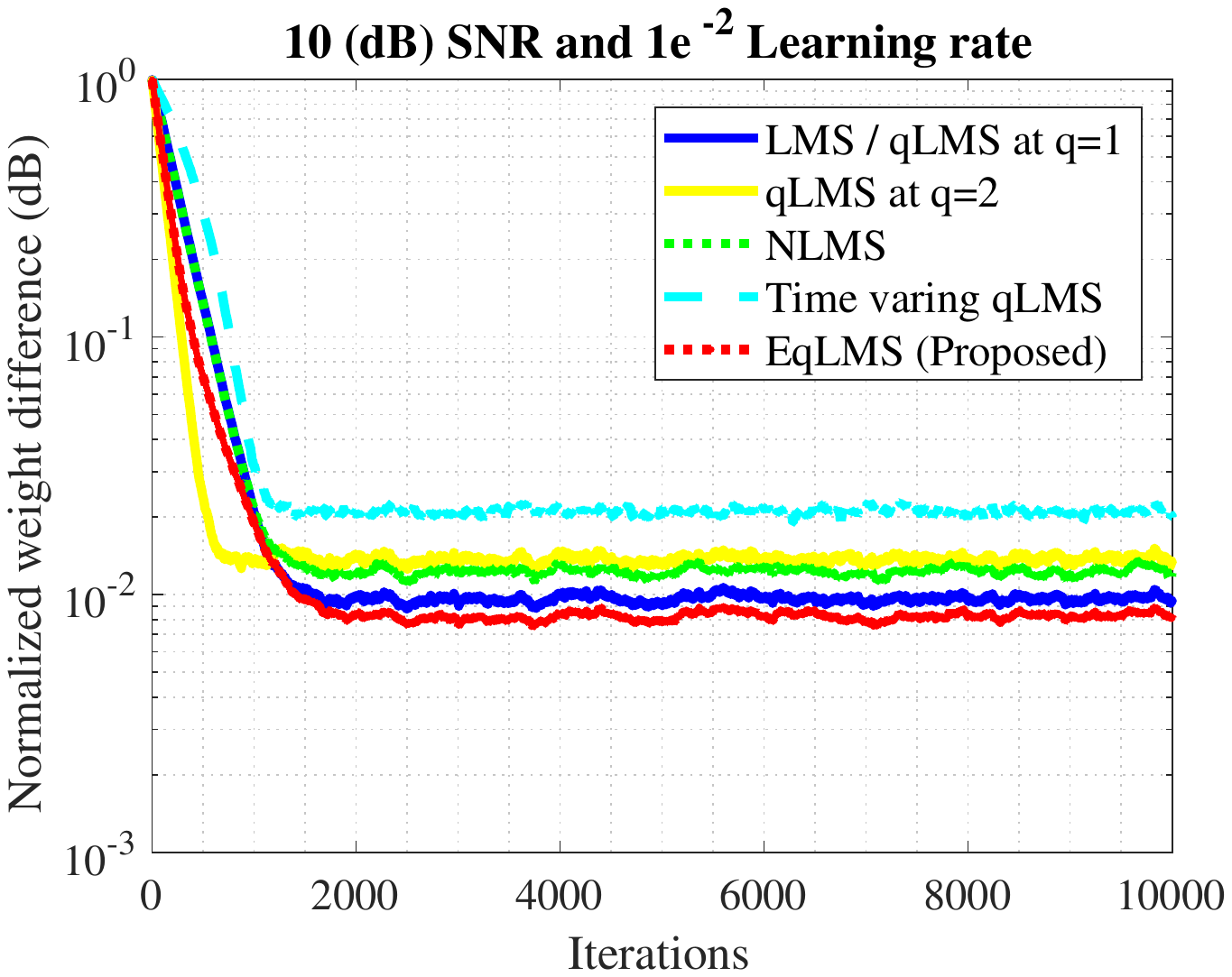}}}
		& 

		\raisebox{-\totalheight}{\centering \fbox{\includegraphics*[scale=0.27,bb=100 240 500 550]{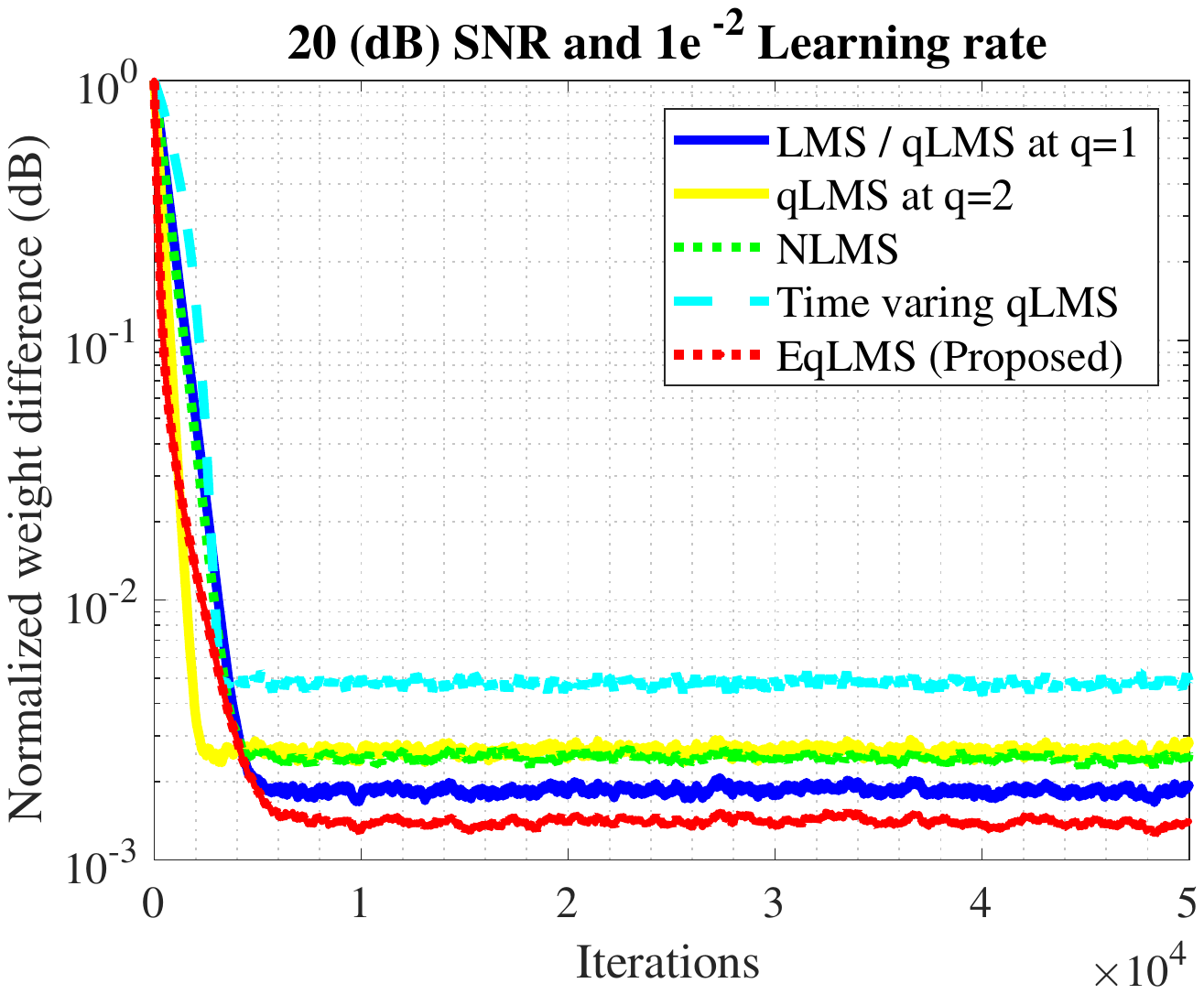}}}
		& 
		\raisebox{-\totalheight}{\centering \fbox{\includegraphics*[scale=0.27,bb=100 240 500 550]{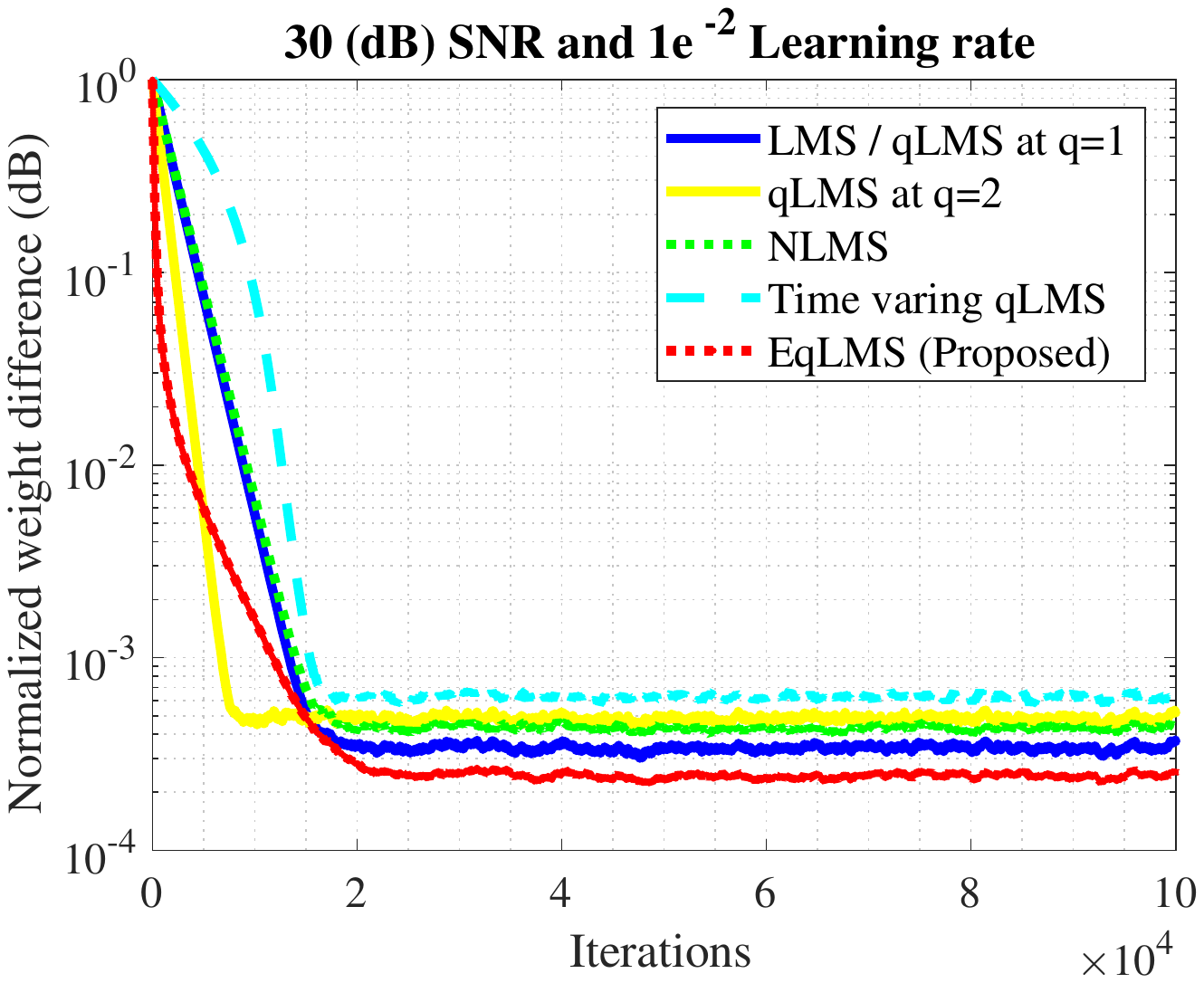}}}
		\\
		d & e & f \\
		\raisebox{-\totalheight}{\centering \fbox{\includegraphics*[scale=0.27,bb=100 240 500 550]{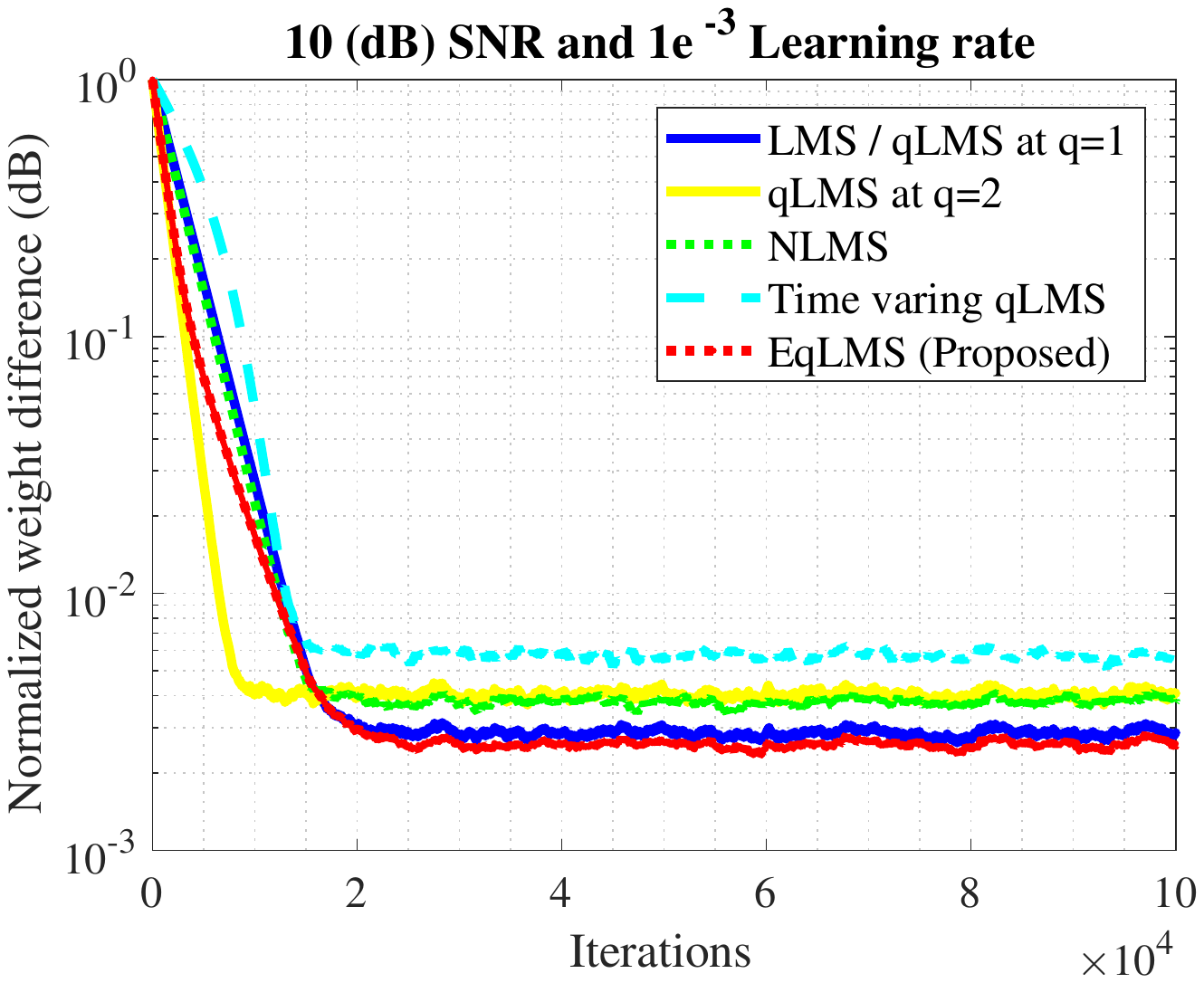}}}
		& 
		\raisebox{-\totalheight}{\centering \fbox{\includegraphics*[scale=0.27,bb=100 240 500 550]{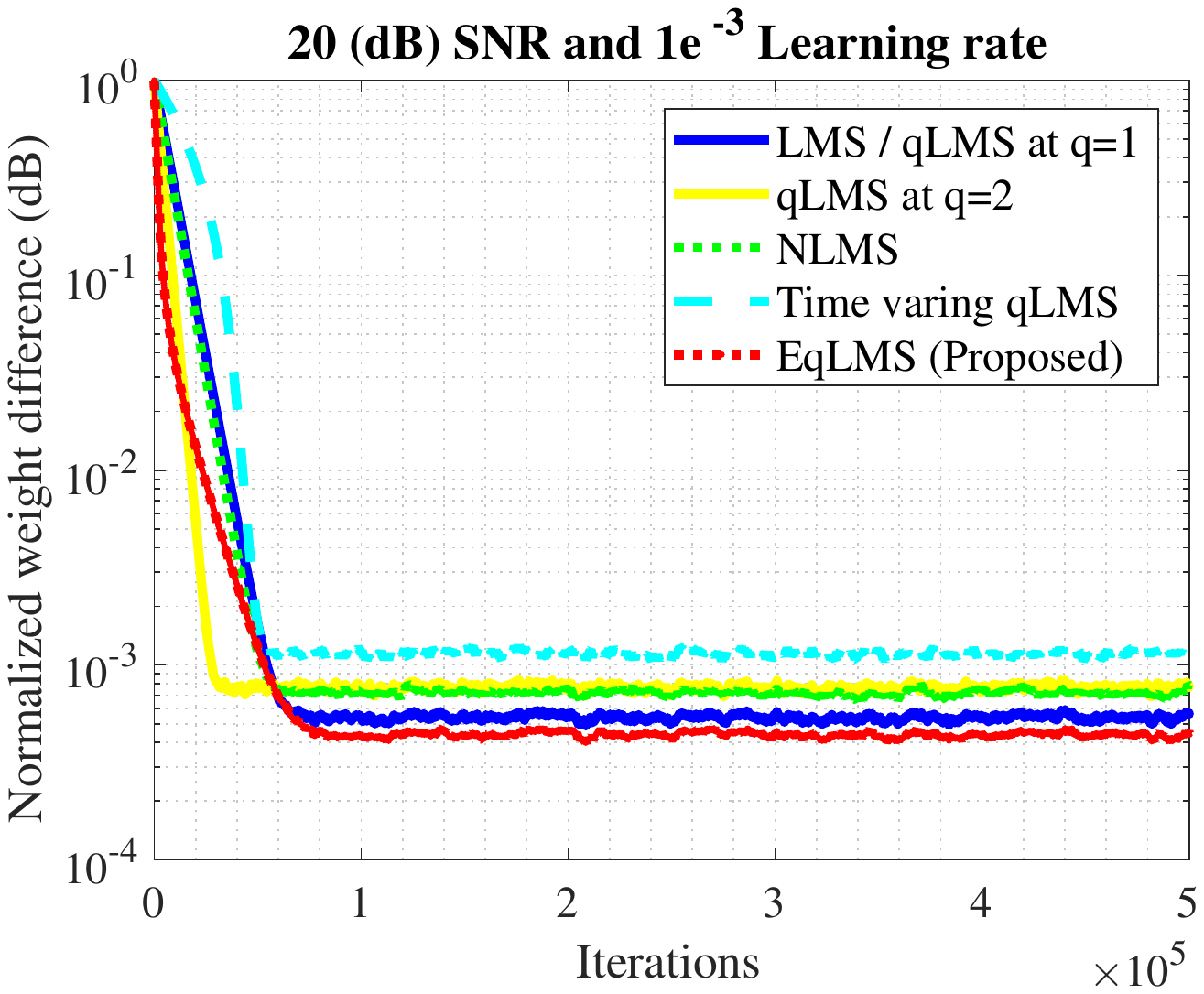}}}
		& 
		\raisebox{-\totalheight}{\centering \fbox{\includegraphics*[scale=0.27,bb=100 240 500 550]{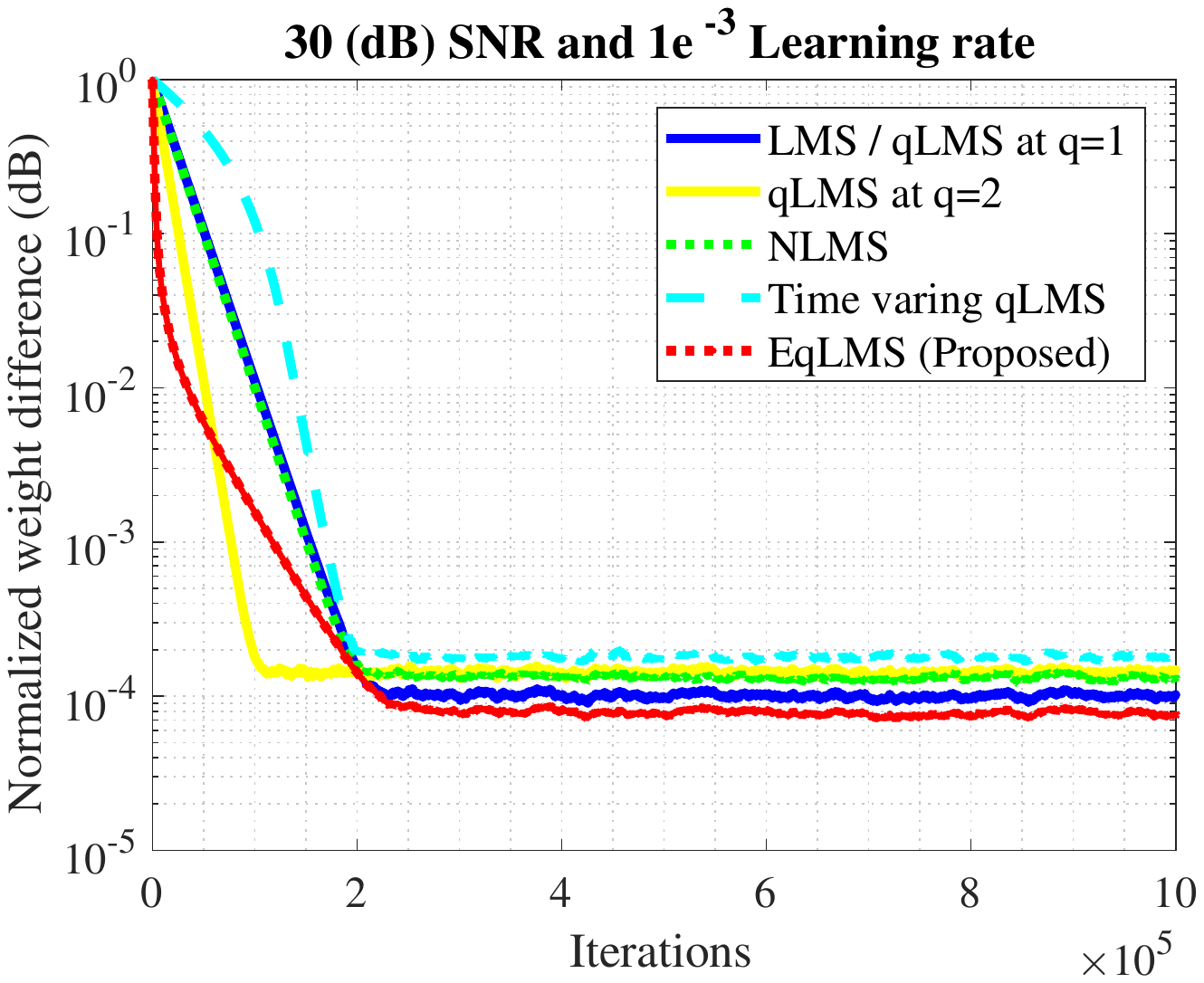}}}			
		\\
		g & h & i 
		\\
	\end{tabular}
	\caption{NWD curves for the LMS/$q$-LMS at $q=1$, $q$-LMS at $q=2$, time-varying $q$-LMS, NLMS, and the $Eq$-LMS. normalized weight deviation with learning rate and SNR of (a) $1e^{-1}$, $10$ dB, (b) $1e^{-1}$, $20$ dB,(c) $1e^{-1}$, $30$ dB, (d) $1e^{-2}$, $10$ dB, (e) $1e^{-2}$, $20$ dB, (f) $1e^{-2}$, $30$ dB, (g) $1e^{-3}$, $10$ dB, (h) $1e^{-3}$, $20$ dB, and (i) $1e^{-3}$, $30$ dB.}
	\label{convergence_scenario}
\end{figure}

The relevant NWD curves with three different SNR values are delineated in Fig. \ref{convergence_scenario}.  From the Fig. \ref{convergence_scenario} (g), (h), and (i), it can be seen that the proposed $Eq$-LMS produced the best performance under all three conditions: (1) for the SNR value of $10$ dB, it outperformed the LMS/$q$-LMS at $q=1$, $q$-LMS at $q=2$, Time-varying $q$-LMS, and the NLMS by the NWD value of $-0.43$ dB, $-1.92$ dB, $-3.45$ dB, and $-1.67$ dB, respectively, (2) for the SNR value of $20$ dB, it surpassed the above-mentioned algorithms by the NWD value of $-0.9$ dB, $-2.43$ dB, $-4.2$ dB, and $-2.15$ dB, respectively, and (3) for the SNR value of $30$ dB, the above mentioned algorithms were outperformed by a margin of $-1.13$ dB, $-2.66$ dB, $-3.68$, and $-2.33$ dB, respectively.  Note that the proposed $Eq$-LMS algorithm showed the lowest steady state error in all conditions while the $q$-LMS at $q=2$ showed faster convergence with greater steady state error than the proposed $Eq$-LMS.  With the above discussed settings, results for the channel estimation problem are summarized in Table \ref{Res_tab_R_scenario3}.  

\begin{table}[H]
	\centering
	\caption{Evaluation protocol 3: Results of various approaches for an equal convergence rate.}
	\label{Res_tab_R_scenario3}
	\begin{tabular}{|l|c|c|c|}
		\hline
		\multicolumn{1}{|c|}{\multirow{2}{*}{\textbf{Algorithm}}} & \multicolumn{3}{c|}{\textbf{Steady-state NWD (dB)}} \\ \cline{2-4} 
		\multicolumn{1}{|c|}{} & \textbf{10 dB SNR} & \textbf{20 dB SNR} & \textbf{30 dB SNR} \\ \hline
		\textbf{LMS/$q$-LMS ($q=1$)} & -20.16 & -27.34 & -34.67 \\ \hline
		\textbf{$q$-LMS ($q=2$)} & -18.92 & -26.13 & -33.45 \\ \hline		
		\textbf{Time-varying $q$-LMS} & -16.76 & -23.22 & -32.04 \\ \hline
		\textbf{Normalized LMS} & -19.10 & -26.11 & -33.64 \\ \hline
		\textbf{Proposed $Eq$-LMS} & \textbf{-20.84} & \textbf{-28.56} & \textbf{-36.11} \\ \hline
	\end{tabular}
\end{table}

\subsection{Experiments to Evaluate the Convergence Performance}
\subsubsection{Evaluation Protocol 1: Fast Convergence}
For the comparison of convergence performances,  all algorithms were setup for equal steady-state error.  The learning rate (step size) configurations for equal steady-state (Evaluation Protocol 1) is shown in Table (\ref{con_tab_S_scenario1}).  Learning rate for proposed $Eq$-LMS has been set according to three evaluation protocols.
\begin{table}[H]
	\centering
\caption{Evaluation protocol 1: Configuration of learning rates of different approaches for an equal steady-state error.}
\label{con_tab_S_scenario1}
	\begin{tabular}{|l|c|c|c|}
		\hline
		\multicolumn{1}{|c|}{\multirow{2}{*}{\textbf{Algorithm}}} & \multicolumn{3}{c|}{\textbf{Learning Rate $\mu$}} \\ \cline{2-4} 
		\multicolumn{1}{|c|}{} & \textbf{10 dB SNR} & \textbf{20 dB SNR} & \textbf{30 dB SNR} \\ \hline
		\textbf{LMS/$q$-LMS} & $3\times10^{-2}$ & $8.9\times10^{-3}$ & $2.7\times10^{-3}$ \\ \hline
		\textbf{Time-varying $q$-LMS} & $3.3\times10^{-2}$ & $8.8\times10^{-3}$ & $3.1\times10^{-3}$ \\ \hline
		\textbf{Normalized LMS} & $1\times10^{-1}$ & $2.8\times10^{-2}$ & $8.5\times10^{-3}$ \\ \hline
		\textbf{Proposed $Eq$-LMS} & $1\times10^{-1}$ & $1\times10^{-1}$ & $1\times10^{-1}$ \\ \hline
	\end{tabular}
\end{table}

The relevant normalized weight difference (NWD) curves with three different SNR values are depicted in Fig. \ref{steady_scenario}.  From the Fig. \ref{steady_scenario} (a), (b), and (c), it can be seen that the proposed $Eq$-LMS algorithm produced the best results under all three conditions: (1) for the SNR value of $10$ dB, algorithms are run for $1000$ iterations.  The convergence point of the proposed $Eq$-LMS is reached at $120^{th}$ iteration, $q$-LMS at ($q=2$) converged on the $80^{th}$ iteration but its steady state error is much larger compared to the proposed $Eq$-LMS, (2) for the SNR value of $20$ dB, algorithms are run for $5000$ iterations.  Note that the proposed $Eq$-LMS algorithm outperformed all competing approaches by converging in only $400$ iterations. The $q$-LMS ($q=2$) was unable to reach the given error-floor and took $400$ iterations to reach a much higher error, and (3) for the SNR value of $30$ dB, algorithms are run for $10000$ iterations, the proposed $Eq$-LMS algorithm took the least number of iterations by converging at the $1600^{th}$ iteration.  The proposed EqLMS shows best performance in terms of steady state error and convergence rate.  Thus, showing the best overall performance.  With the above discussed settings, results for the channel estimation problem are summarized in Table \ref{Res_tab_S_scenario1}.  

\begin{table}[H]
	\centering
	\caption{Evaluation protocol 1: Results of various approaches for an equal steady-state error.}
\label{Res_tab_S_scenario1}
	\begin{tabular}{|l|c|c|c|}
		\hline
		\multicolumn{1}{|c|}{\multirow{2}{*}{\textbf{Algorithm}}} & \multicolumn{3}{c|}{\textbf{Convergence point (number of Iterations $\times1000$))}} \\ \cline{2-4} 
		\multicolumn{1}{|c|}{} & \textbf{10 dB SNR} & \textbf{20 dB SNR} & \textbf{30 dB SNR} \\ \hline
		\textbf{LMS/$q$-LMS at ($q=1$)} & $0.20$ & $0.70$ & $2.70$ \\ \hline
		\textbf{$q$-LMS at ($q=2$)} & $0.08$ & $0.40$ & $1.55$ \\ \hline
		\textbf{Time-varying $q$-LMS} & $0.64$ & $2.80$ & $7.20$ \\ \hline
		\textbf{Normalized LMS} & $0.23$ & $1.10$ & $4.40$ \\ \hline
		\textbf{Proposed $Eq$-LMS} & $0.12$ & $0.40$ & $1.60$ \\ \hline
	\end{tabular}
\end{table}

\subsubsection{Evaluation protocol 2: Medium Convergence}
The learning rate (step size) configuration for equal steady-state (Evaluation Protocol 2) is shown in configuration Table (\ref{con_tab_S_scenario2}).  Learning rate for proposed $Eq$-LMS has been set according to three evaluation protocols.
\begin{table}[H]
	\centering
	\caption{Evaluation protocol 2: Configuration of learning rates of different approaches for an equal steady-state error.}
\label{con_tab_S_scenario2}
	\begin{tabular}{|l|c|c|c|}
		\hline
		\multicolumn{1}{|c|}{\multirow{2}{*}{\textbf{Algorithm}}} & \multicolumn{3}{c|}{\textbf{Learning Rate $\mu$}} \\ \cline{2-4} 
		\multicolumn{1}{|c|}{} & \textbf{10 dB SNR} & \textbf{20 dB SNR} & \textbf{30 dB SNR} \\ \hline
		\textbf{LMS/$q$-LMS} & $3\times10^{-3}$ & $8.8\times10^{-4}$ & $2.7\times10^{-4}$ \\ \hline
		\textbf{Time-varying $q$-LMS} & $3.3\times10^{-3}$ & $9.3\times10^{-4}$ & $3.1\times10^{-4}$ \\ \hline
		\textbf{Normalized LMS} & $9\times10^{-2}$ & $2.72\times10^{-3}$ & $8.5\times10^{-4}$ \\ \hline
		\textbf{Proposed $Eq$-LMS} & $1\times10^{-2}$ & $1\times10^{-2}$ & $1\times10^{-2}$ \\ \hline
	\end{tabular}
\end{table}

The relevant NWD curves with three different SNR values are depicted in Fig. \ref{steady_scenario}.  From the Fig. \ref{steady_scenario} (d), (e), and (f), it can be seen that the proposed $Eq$-LMS algorithm produced the best results under all three conditions: (1) for the SNR value of $10$ dB, algorithms are run for $10000$ iterations, the convergence point of the proposed $Eq$-LMS is reached at $1500^{th}$ iteration, the $q$-LMS at ($q=2$) converged on the $1100^{th}$ iteration but its steady state error is much larger than the proposed $Eq$-LMS, (2) for the SNR value of $20$ dB, algorithms are run for $50000$ iterations, the proposed $Eq$-LMS algorithm outperformed all competing approaches in terms of convergence point with least steady state error, and (3) for the SNR value of $30$ dB, algorithms are run for $100000$ iterations, the proposed $Eq$-LMS algorithm converged on $19000^{th}$ iteration, it showed best performance in terms of steady state error and convergence rate.  Thus, showing the best overall performance.  With the above discussed settings, results for the channel estimation problem are summarized in Table \ref{Res_tab_S_scenario2}. 
\begin{table}[H]
	\centering
\caption{Evaluation protocol 2: Results of various approaches for an equal steady-state error.}
\label{Res_tab_S_scenario2}
	\begin{tabular}{|l|c|c|c|}
		\hline
		\multicolumn{1}{|c|}{\multirow{2}{*}{\textbf{Algorithm}}} & \multicolumn{3}{c|}{\textbf{Convergence point (number of Iterations $\times1000$))}} \\ \cline{2-4} 
		\multicolumn{1}{|c|}{} & \textbf{10 dB SNR} & \textbf{20 dB SNR} & \textbf{30 dB SNR} \\ \hline
		\textbf{LMS/$q$-LMS at ($q=1$)} & $2$ & $9$ & $32$ \\ \hline
		\textbf{$q$-LMS at ($q=2$)} & $1.1$ & $5.9$ & $19$ \\ \hline
		\textbf{Time-varying $q$-LMS} & $6.2$ & $28$ & $80$ \\ \hline
		\textbf{Normalized LMS} & $3.1$ & $14$ & $52$ \\ \hline
		\textbf{Proposed $Eq$-LMS} & $1.5$ & $6$ & $19$ \\ \hline
	\end{tabular}
\end{table}

\subsubsection{Evaluation protocol 3: Slow Convergence}
The learning rate (step size) configuration for equal steady-state (Evaluation protocol 3) is shown in configuration Table (\ref{con_tab_S_scenario3}).  Learning rate for the proposed $Eq$-LMS has been set according to three evaluation protocols.
\begin{table}[H]
	\centering
\caption{Evaluation protocol 3: Configuration of learning rates of different approaches for an equal steady-state error.}
\label{con_tab_S_scenario3}
	\begin{tabular}{|l|c|c|c|}
		\hline
		\multicolumn{1}{|c|}{\multirow{2}{*}{\textbf{Algorithm}}} & \multicolumn{3}{c|}{\textbf{Learning Rate $\mu$}} \\ \cline{2-4} 
		\multicolumn{1}{|c|}{} & \textbf{10 dB SNR} & \textbf{20 dB SNR} & \textbf{30 dB SNR} \\ \hline
		\textbf{LMS/$q$-LMS} & $3\times10^{-4}$ & $8.9\times10^{-4}$ & $2.7\times10^{-5}$ \\ \hline
		\textbf{Time-varying $q$-LMS} & $3.3\times10^{-4}$ & $9\times10^{-5}$ & $3.2\times10^{-5}$ \\ \hline
		\textbf{Normalized LMS} & $1\times10^{-3}$ & $2.8\times10^{-4}$ & $8.5\times10^{-5}$ \\ \hline
		\textbf{Proposed $Eq$-LMS} & $1\times10^{-3}$ & $1\times10^{-3}$ & $1\times10^{-3}$ \\ \hline
	\end{tabular}
\end{table}

\begin{figure}[]
	
	\begin{tabular}{>{\centering}m{3.7cm} >{\centering}m{3.7cm} >{\centering\arraybackslash}m{3.7cm}}
		
		\raisebox{-\totalheight}{\centering \fbox{\includegraphics*[scale=0.27,bb=100 240 500 550]{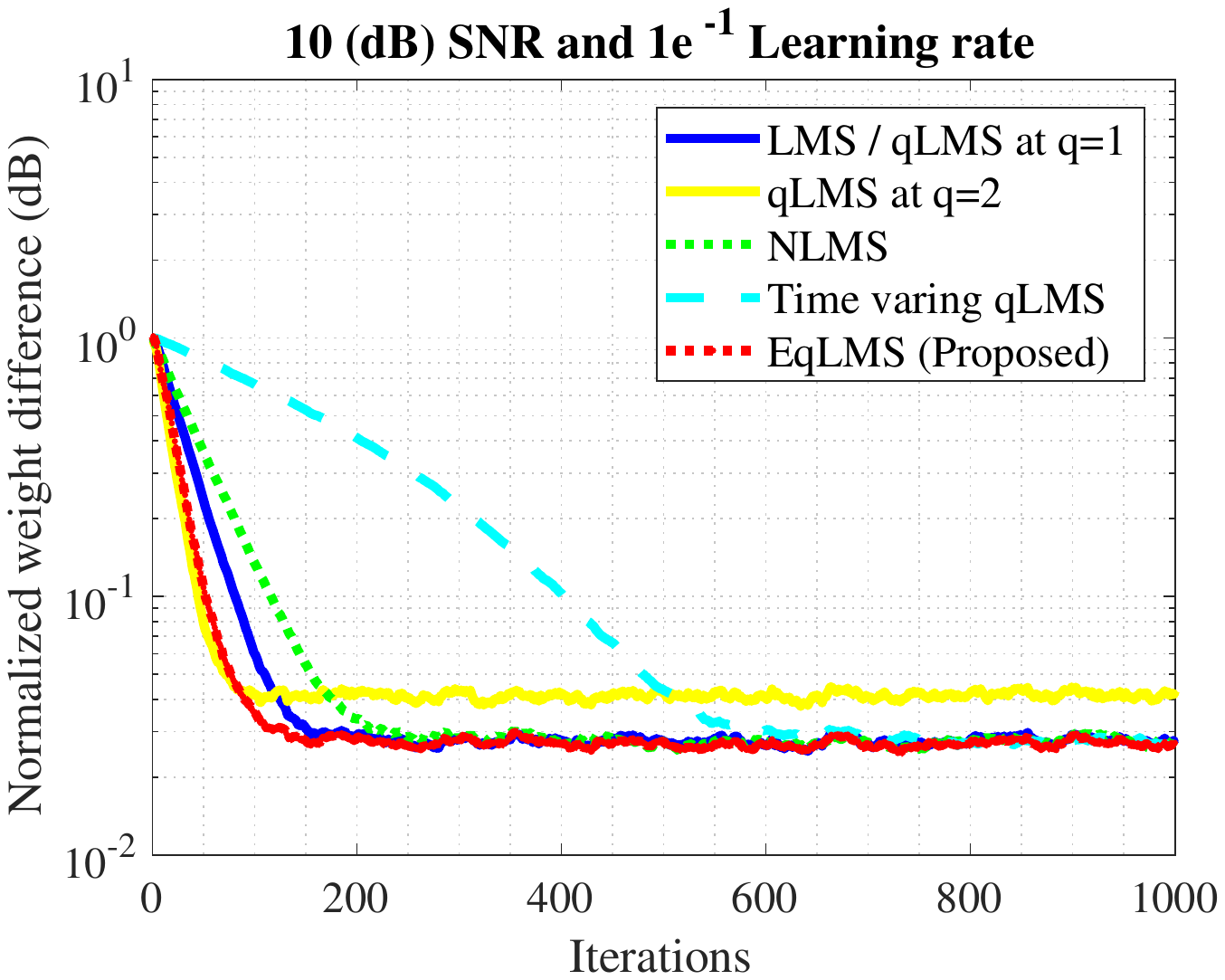}}}
		& 
		
		\raisebox{-\totalheight}{\centering \fbox{\includegraphics*[scale=0.27,bb=100 240 500 550]{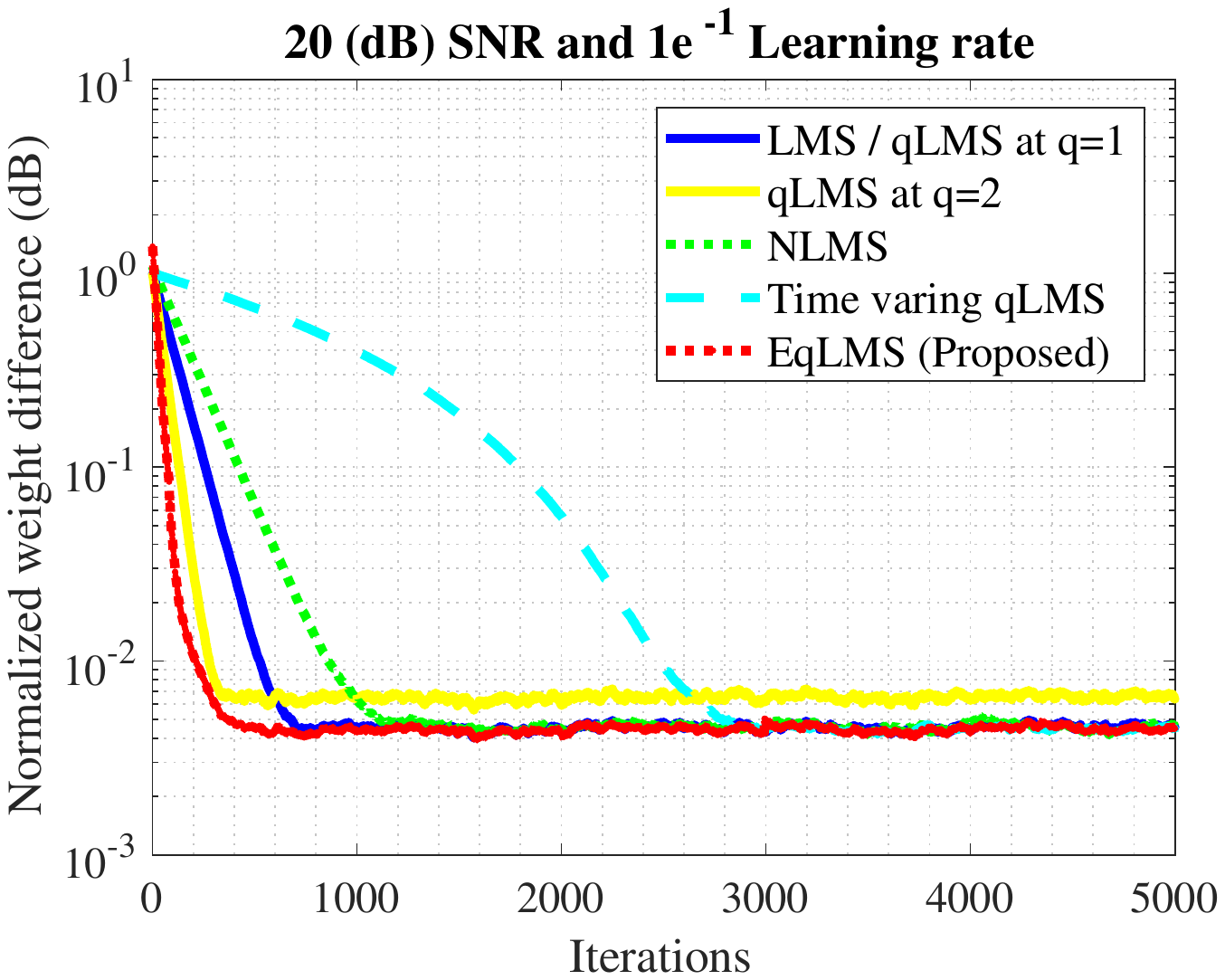}}}
		& 
		
		\raisebox{-\totalheight}{\centering \fbox{\includegraphics*[scale=0.27,bb=100 240 500 550]{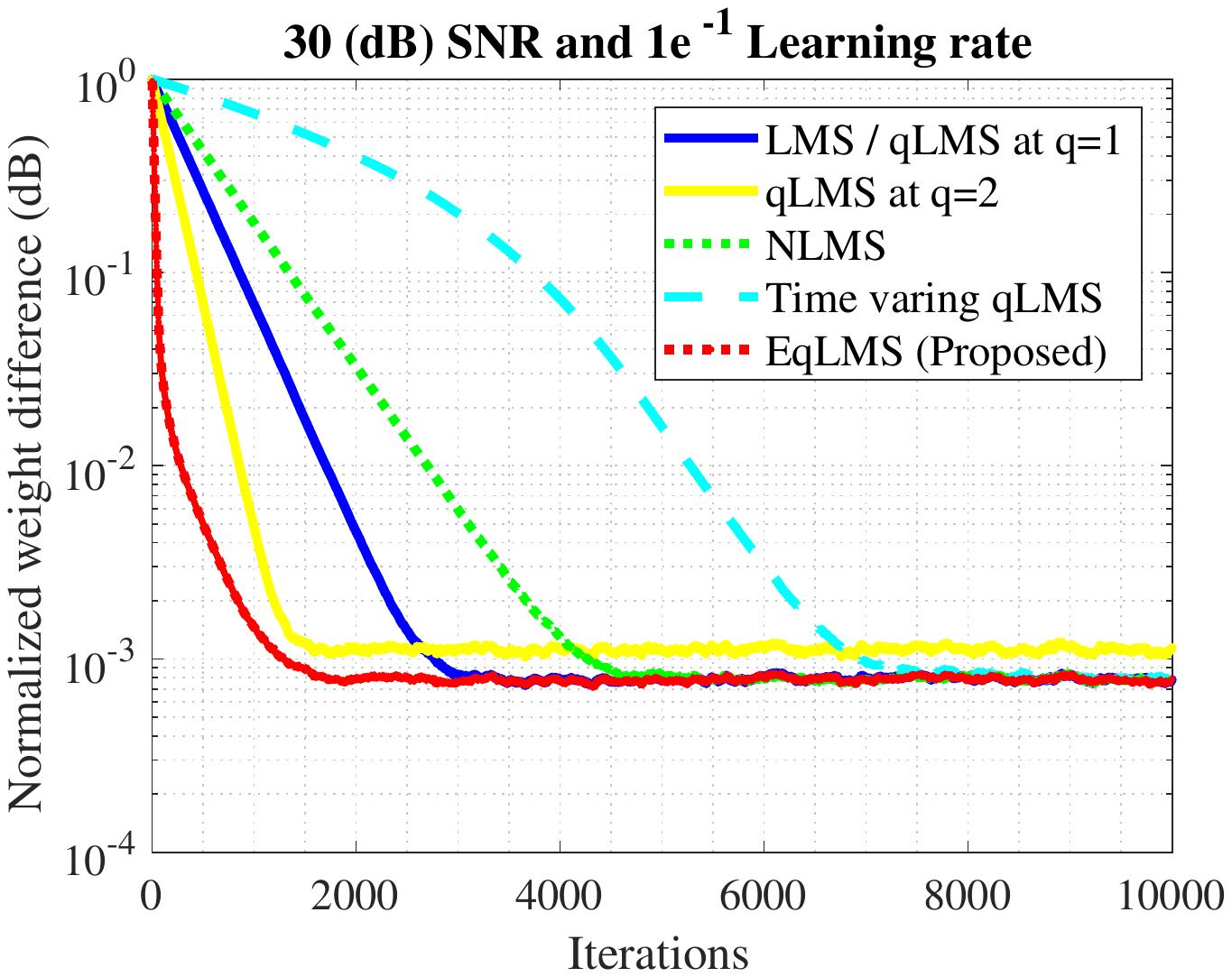}}}
		\\ 
		a & b & c \\
		
		\raisebox{-\totalheight}{\centering \fbox{\includegraphics*[scale=0.27,bb=100 240 500 550]{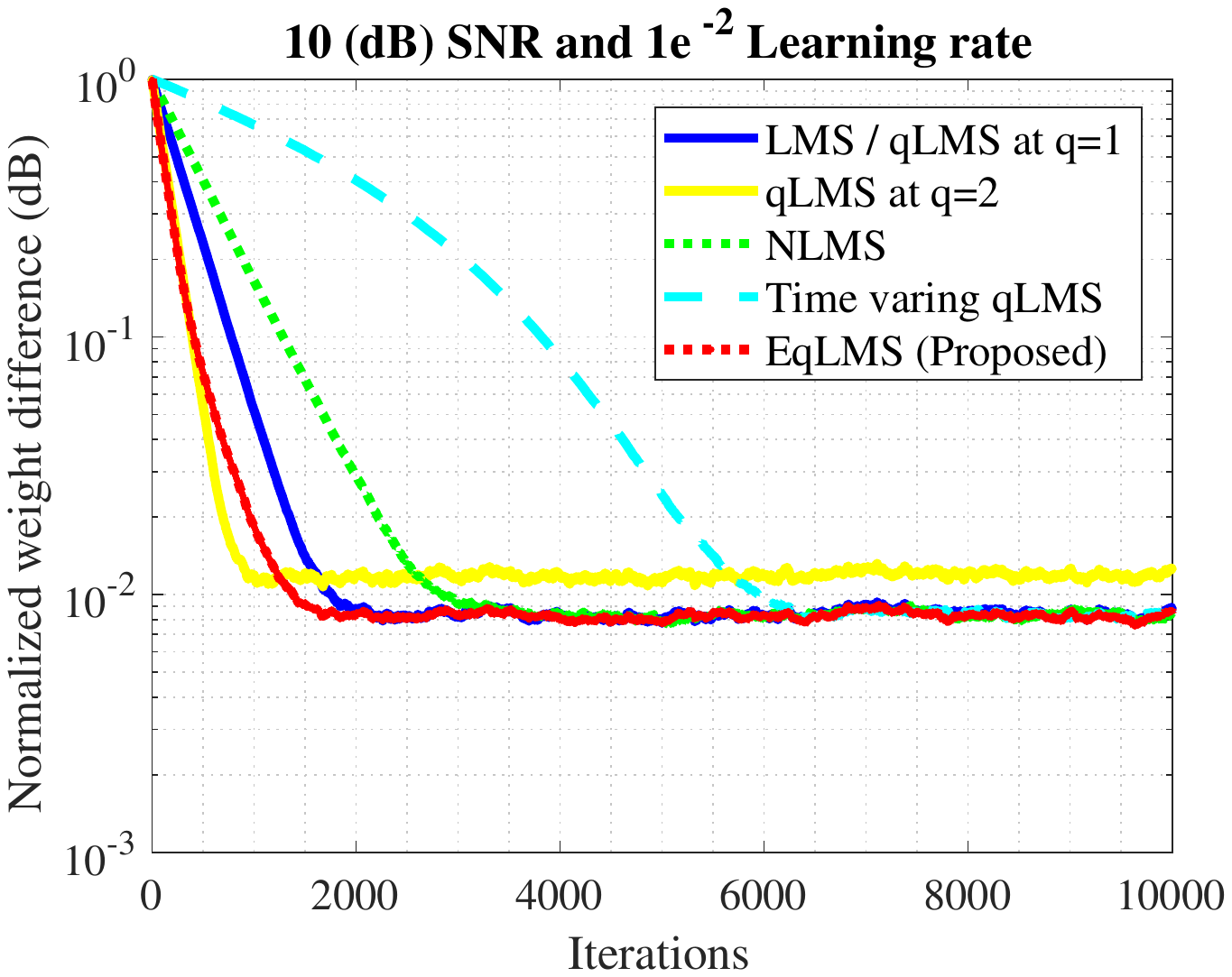}}}
		& 
		
		\raisebox{-\totalheight}{\centering \fbox{\includegraphics*[scale=0.27,bb=100 240 500 550]{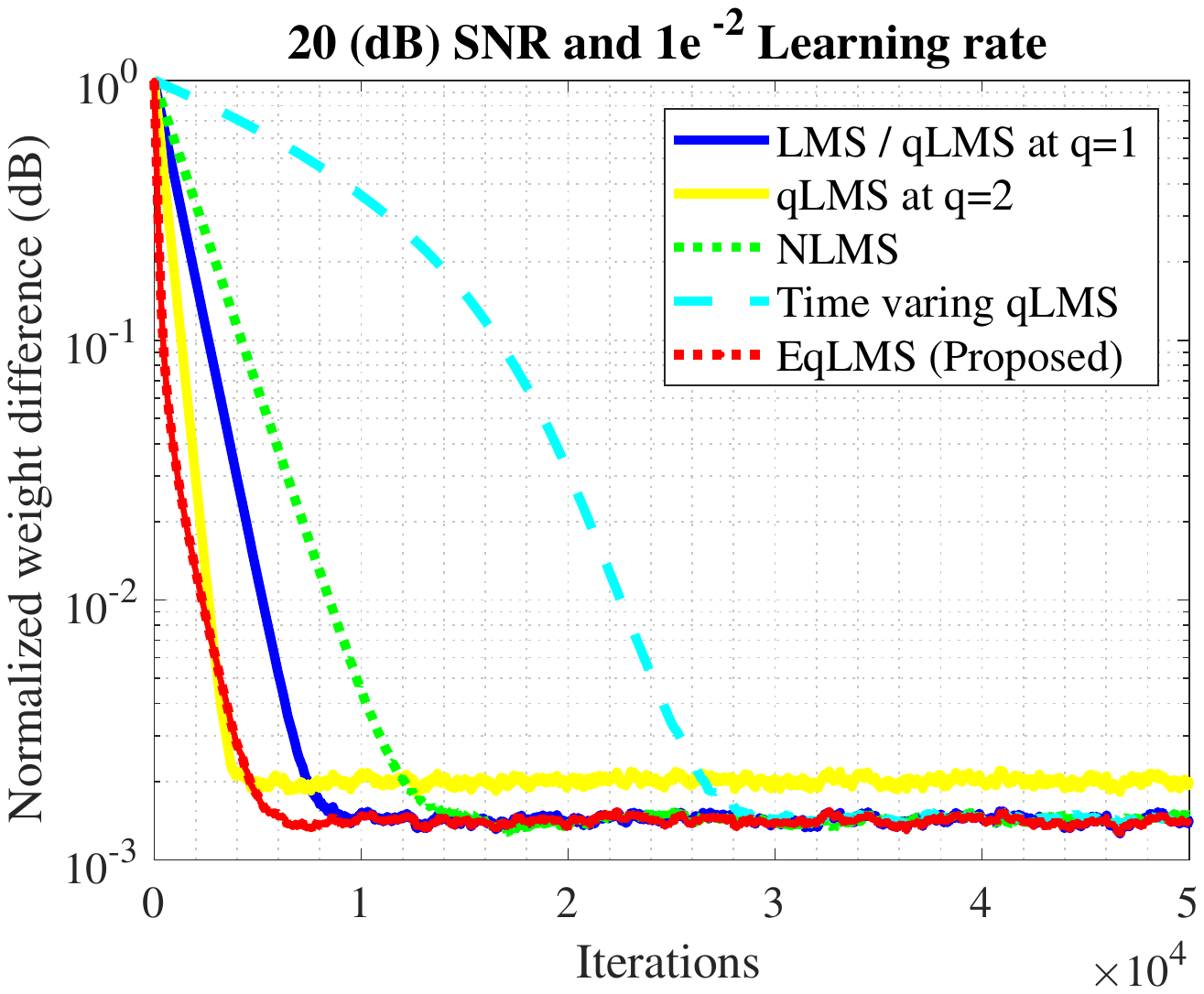}}}
		& 
		
		\raisebox{-\totalheight}{\centering \fbox{\includegraphics*[scale=0.27,bb=100 240 500 550]{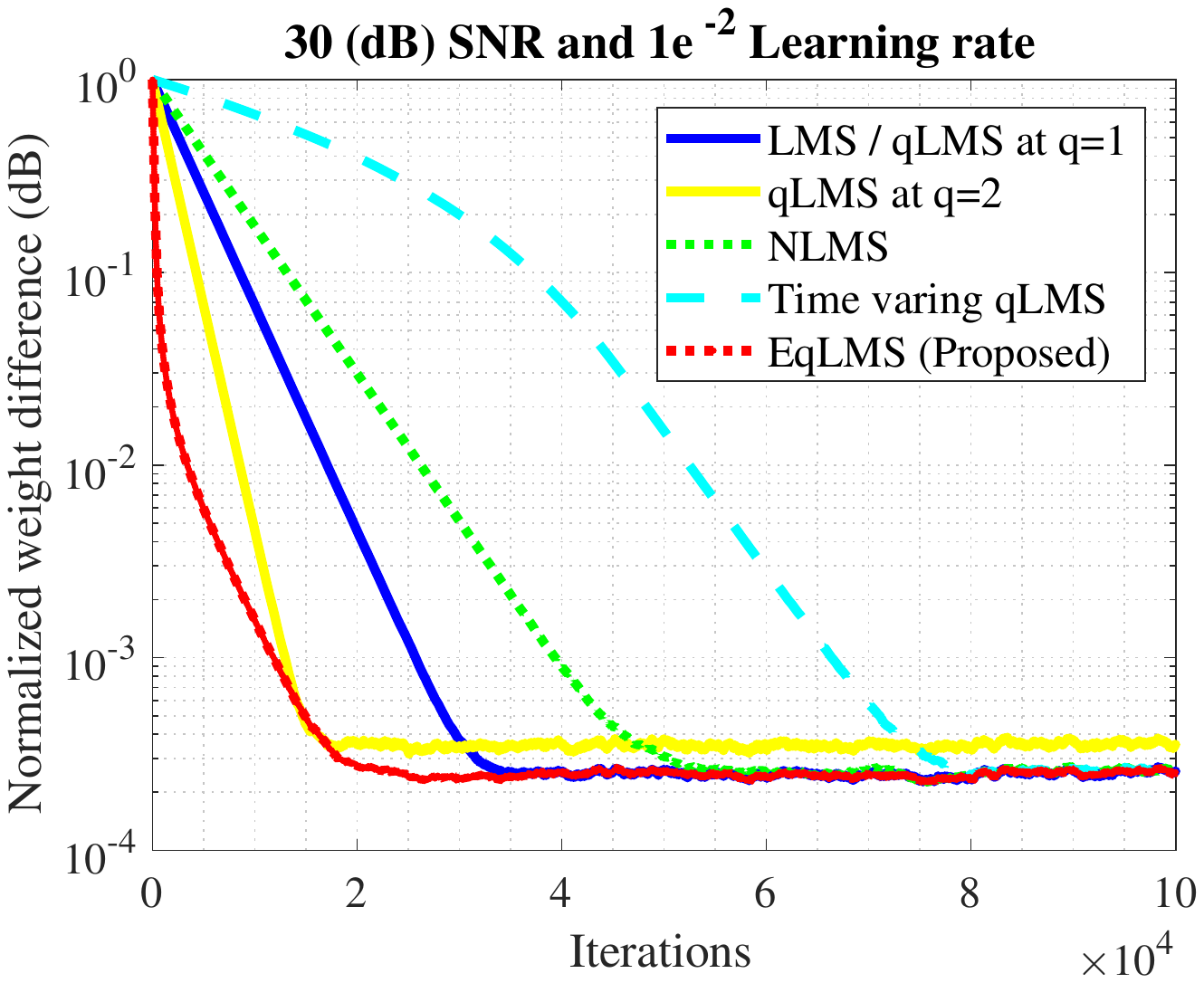}}}
		\\
		d & e & f \\
		
		\raisebox{-\totalheight}{\centering \fbox{\includegraphics*[scale=0.27,bb=100 240 500 550]{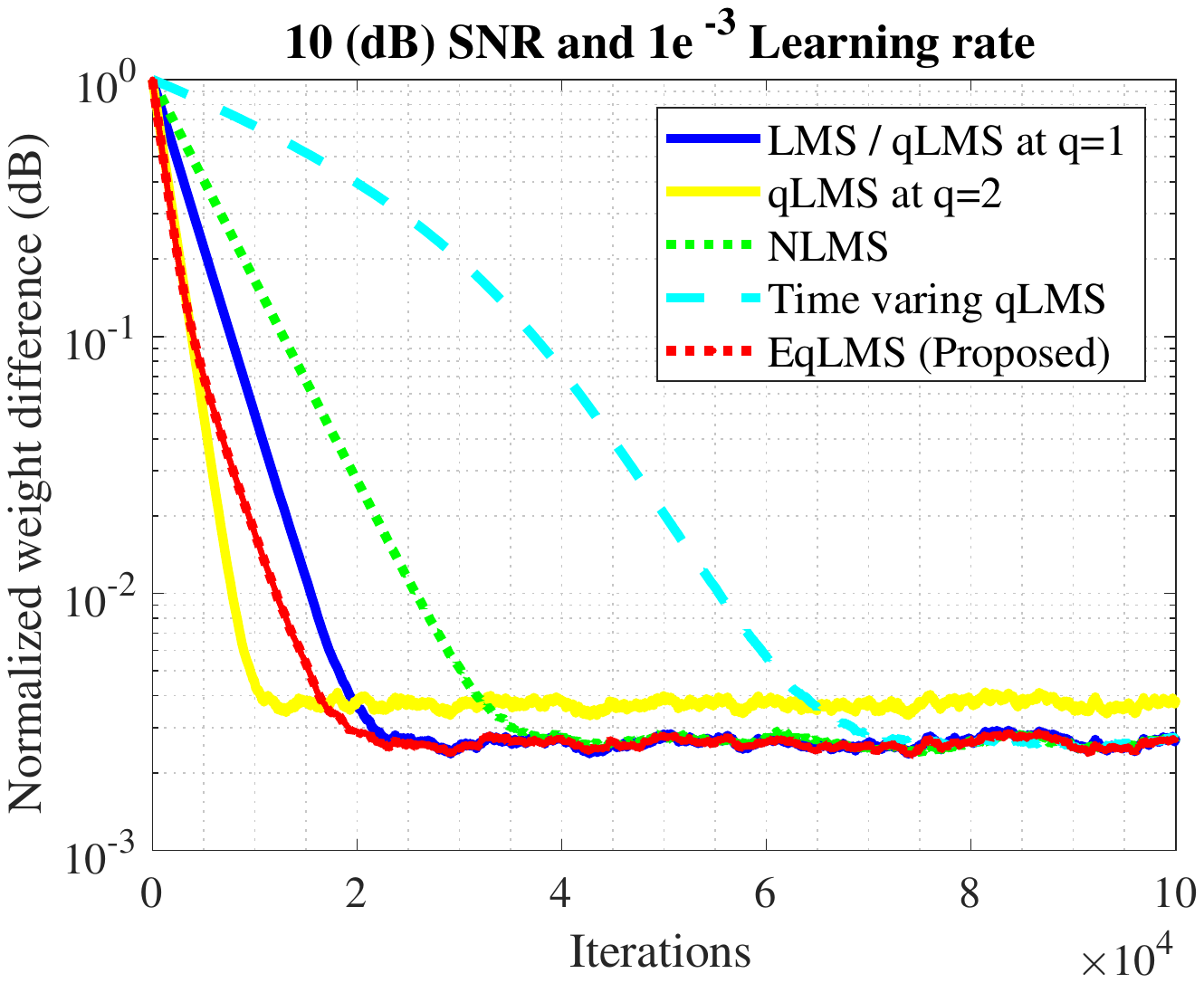}}}
		& 
		
		\raisebox{-\totalheight}{\centering \fbox{\includegraphics*[scale=0.27,bb=100 240 500 550]{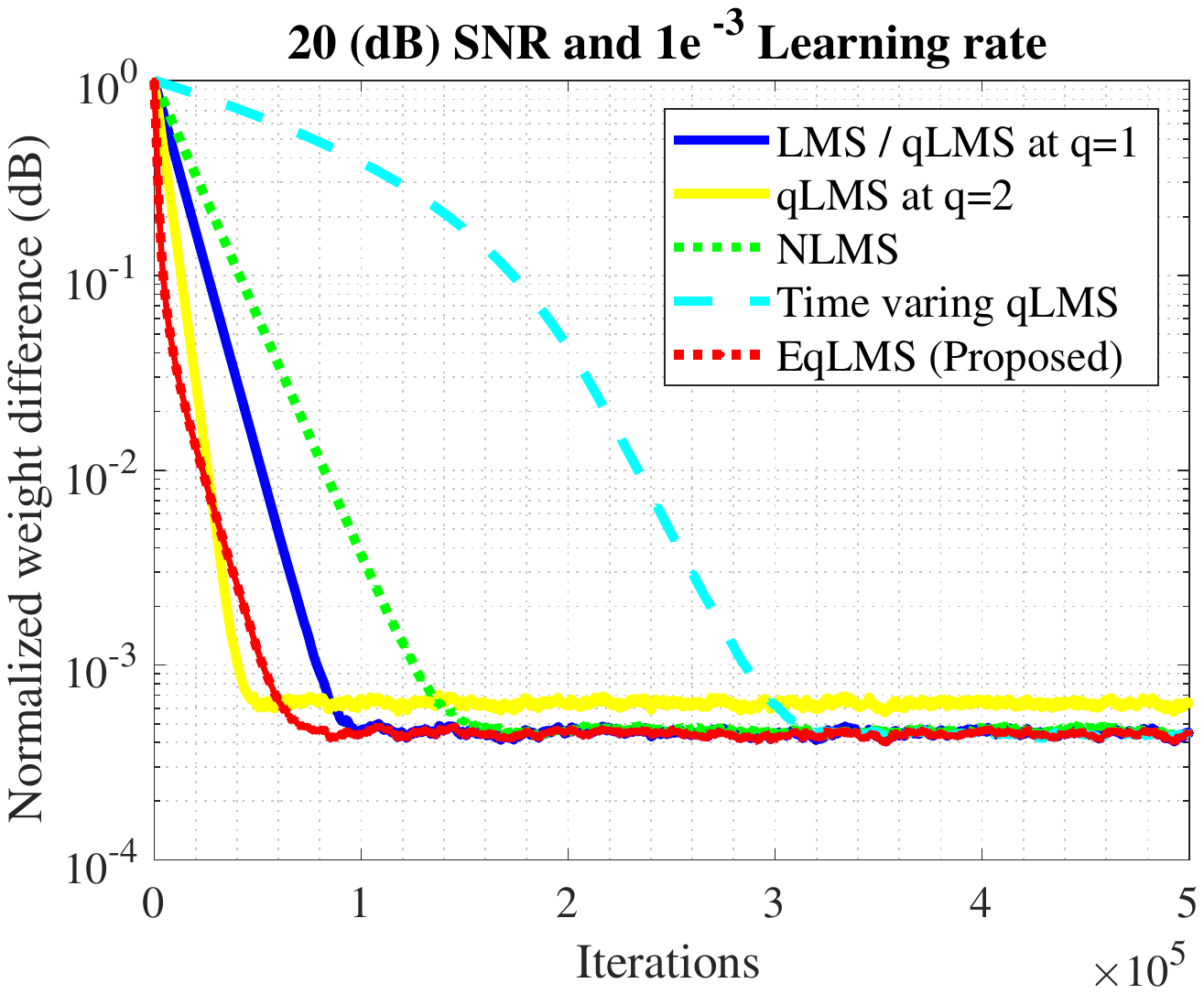}}}
		& 
		
		\raisebox{-\totalheight}{\centering \fbox{\includegraphics*[scale=0.27,bb=100 240 500 550]{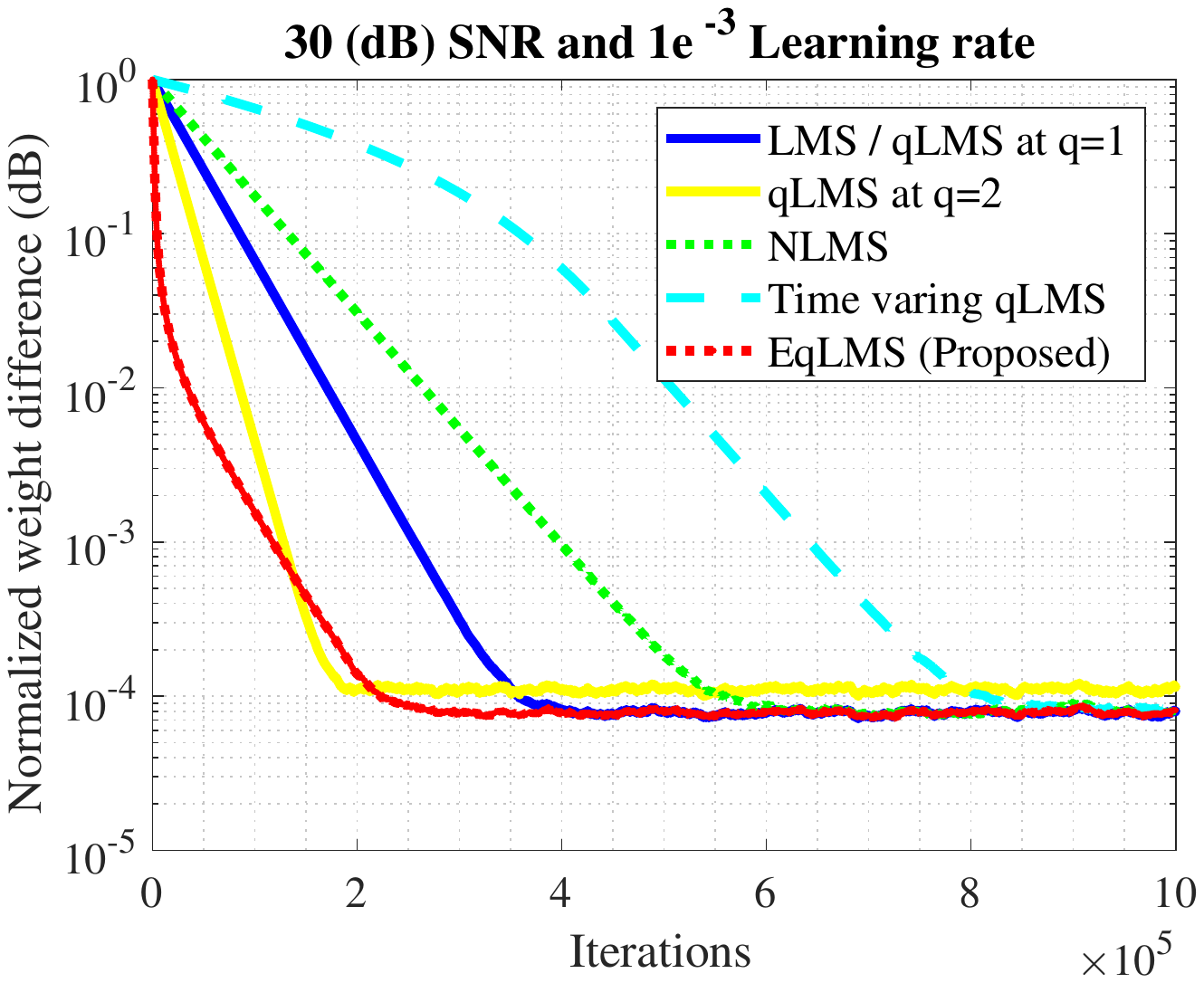}}}
		\\
		g & h & i 
		\\
	\end{tabular}
	\caption{NWD curves for the LMS/$q$-LMS at $q=1$, $q$-LMS at $q=2$, time-varying $q$-LMS, NLMS, and the $Eq$-LMS. normalized weight deviation with learning rate and SNR of (a) $1e^{-1}$, $10$ dB, (b) $1e^{-1}$, $20$ dB,(c) $1e^{-1}$, $30$ dB, (d) $1e^{-2}$, $10$ dB, (e) $1e^{-2}$, $20$ dB, (f) $1e^{-2}$, $30$ dB, (g) $1e^{-3}$, $10$ dB, (h) $1e^{-3}$, $20$ dB, and (i) $1e^{-3}$, $30$ dB.}
	\label{steady_scenario}
	
\end{figure}

The relevant NWD curves with three different SNR values are depicted in Fig. \ref{steady_scenario}.  From the Fig. \ref{steady_scenario} (g), (h), and (i), it can be seen that the proposed $Eq$-LMS algorithm produced the best results under all three conditions: (1) for the SNR value of $10$ dB, algorithms are run for $100000$ iterations, the convergence point of the proposed $Eq$-LMS is reached at $21000^{th}$ iteration, the $q$-LMS at ($q=2$), converged at the $12000^{th}$ iteration but its steady state error is much larger than the proposed $Eq$-LMS, (2) for the SNR value of $20$ dB, the algorithms are run for $500000$ iterations, the proposed $Eq$-LMS algorithm outperformed all competing approaches in terms of convergence point with least steady state error, and (3) for the SNR value of $30$ dB, the algorithms are run for $1000000$ iterations, the proposed $Eq$-LMS algorithm converged on $200000^{th}$ iteration.  The proposed $Eq$-LMS showed the best performance in terms of steady state error and convergence rate.   With the above discussed settings, results for the channel estimation problem are summarized in Table \ref{Res_tab_S_scenario3}.   
\begin{table}[H]
	\centering
	\caption{Evaluation protocol 3: Results of various approaches for an equal steady-state error.}
\label{Res_tab_S_scenario3}
	\begin{tabular}{|l|c|c|c|}
		\hline
		\multicolumn{1}{|c|}{\multirow{2}{*}{\textbf{Algorithm}}} & \multicolumn{3}{c|}{\textbf{Convergence point (number of Iterations $\times1000$))}} \\ \cline{2-4} 
		\multicolumn{1}{|c|}{} & \textbf{10 dB SNR} & \textbf{20 dB SNR} & \textbf{30 dB SNR} \\ \hline
		\textbf{LMS/qLMS at (q=1)} & $23$ & $100$ & $370$ \\ \hline
		\textbf{qLMS at (q=2)} & $12$ & $60$ & $200$ \\ \hline
		\textbf{Time-varying qLMS} & $72$ & $28$ & $840$ \\ \hline
		\textbf{Normalized LMS} & $36$ & $320$ & $600$ \\ \hline
		\textbf{Proposed Eq-LMS} & $21$ & $70$ & $240$ \\ \hline
	\end{tabular}
\end{table}

\section{Conclusion}\label{Sec:Con}
In this work, we proposed a quantum calculus-based steepest descent algorithm called enhanced $q$-least mean square algorithm ($Eq$-LMS) using a novel concept of error correlation energy.  The proposed algorithm is a parameterless method and unlike the contemporary time varyying $q$-LMS, it does not require additional tuning.  The proposed $Eq$-LMS was compared with the LMS, $q$-LMS, time varying $q$-LMS, and the NLMS algorithms for a problem of linear channel estimation.  Extensive simulation tests were conducted to analyze the convergence and the steady-state performance at three different SNR levels.  For all scenarios, the proposed $Eq$-LMS algorithm comprehensively outperformed the contemporary approaches achieving the best performance in terms of steady-state error and convergence.

%

\end{document}